\documentclass[11pt]{article}
\usepackage{amsfonts}
\usepackage{latexsym,amsmath}
\usepackage{amssymb,array}
\usepackage[symbol]{footmisc}
\usepackage{fancyhdr}
\usepackage[a4paper, portrait, margin=1.1811in]{geometry}
\usepackage{amsthm}
\usepackage{relsize}
\usepackage{graphicx}
\usepackage{booktabs}
\usepackage{lscape}
\usepackage{longtable}
\usepackage{xcolor} 
\usepackage{multirow}
\usepackage{threeparttable}
\usepackage[longnamesfirst,sort]{natbib}
\bibpunct[, ]{(}{)}{;}{a}{,}{,}

\usepackage[colorlinks, citecolor=cyan]{hyperref}
\newtheorem{t1}{Theorem}[section]
\newtheorem{p1}{Proposition}[section]

\begin{document}
	\title{\textbf{Fractional differential entropy and its application in modeling one-dimensional flow velocity}}
	\author{Poulami Paul\and Chanchal Kundu\footnote{\textit{Corresponding
				author e-mail}: \color{cyan}ckundu@rgipt.ac.in;
			\color{cyan}chanchal$_{-}$kundu@yahoo.com.}
		\and Department of Mathematical Sciences\\
		Rajiv Gandhi Institute of Petroleum Technology\\
		Jais 229 304, U.P., India}
	\date{July, 2025}
	\maketitle
	\begin{abstract}
		The fractional order generalization of Shannon entropy proposed by Ubriaco has been studied for discrete distributions. In the current paper, we conduct a detailed study of the continuous analogue of this entropy termed as fractional differential entropy and find some interesting properties which makes it stand out among the existing entropies in literature. The studied entropy measure is evaluated analytically and numerically for some well-known continuous distributions, which will be quite useful in reliability analysis works and other statistical studies of complex systems. Further, it has been used to model the one-dimensional vertical velocity profile of turbulent flows in wide open channels. A one-parametric spatial distribution function is utilized for better estimation of the velocity distribution. The validity of the model has been established using experimental and field data through regression analysis. A comparative study is also presented to show the superiority of the proposed model over the existing entropy-based models. 
	\end{abstract}
	\section{Introduction}
	Intrigued by the thermodynamic concept of entropy as a quantitative measure of disorder of a physical system, Shannon 
	introduced a similar concept in information theory to analyze the randomness or uncertainty associated with information sources and communication systems, enabling the prediction of their output behavior (cf. Shannon, 1948). The uncertainty measure proposed by Shannon for discrete random variables is defined as
	\begin{equation}
		H_S(p) = - \sum_{i=1}^{n} p_i \log p_i,
		\label{E1}
	\end{equation}
	where $ p_i $ is the probability mass function of the discrete random variable. This measure can be directly extended to the continuous case but at the cost of loosing some classical properties.
	If we have an absolutely continuous nonnegative random variable $ X $ with cumulative distribution function (cdf) $ F(\cdot) $ such that $ \frac{dF(x)}{dx} = f(x) $, where $ f(\cdot) $ is the probability density function (pdf), then the Shannon differential entropy of $ X $ is given by 
	\begin{equation}
		H_S(f) = E[-\log f(x)] = - \int_{0}^{\infty} f(x)\log f(x)~dx \label{1.1}  
	\end{equation}where $ E[\cdot] $ gives the expected value of the uncertainty associated with $ f $. 
	Here, $ \log $ stands for the natural logarithm.\\
	\hspace*{0.2in} Thenceforth, the scientific community across the world dived into the work of parametric generalization of the Shannon entropy (\ref{E1}) and (\ref{1.1}) to study the behavior of chaotic physical systems which couldn't be described explicitly by the Shannon entropy. A lot of generalizations have been proposed using the tools of fractional calculus (cf. Machado and Lopes, 2019 and Val$\acute{e}$rio et al., 2013), which allowed the relaxation of some properties of Shannon entropy and thereby making it eligible for describing the complexity and dynamic behavior of physical processes. One may refer to Foroghi et al. (2022), Karci (2016), Lopes and Machado (2020) and references therein for further reading.\\
	\hspace*{0.2in} Ubriaco (2009) generalized the Shannon entropy for discrete random variables and proposed a new definition using the fact that Shannon entropy can be derived by Abe (1998) 
	\begin{equation}
		H_S(p) = \lim\limits_{t\to -1}\frac{d}{dt}\sum\limits_{i} p(i)^{-t}. \label{1.3}	
	\end{equation} Replacing the ordinary derivative $ \frac{d}{dt} $ in (\ref{1.3}) with left Riemann-Liouville fractional derivative $ {}_a^{RL}D_t^{\alpha}, a = -\infty $ and $ 0 \leq \alpha < 1 $ and using the fact that $ p(i)^{-t} $ can be expressed as $ e^{-t\log p(i)} $ within the summation (\ref{1.3}), 
	he came up with a novel fractional order entropy given by 
	\begin{equation}
		H_U^\alpha (p) = \sum_{i=1}^{n}p(i) (-\log p(i))^\alpha;~ \alpha \in [0,1] \label{1.4},
	\end{equation} where, $ {}_{-\infty}^{RL}D_t^{\alpha} f(t) =  \frac{1}{\Gamma (1-\alpha)}\frac{d}{dt}\int\limits_{-\infty}^t \frac{f(u)}{(t-u)^{\alpha}} du,~ 0 \leq \alpha < 1. $  \\
	\hspace*{0.2in} Later on, Machado (2014) generalized the Ubriaco entropy for fractional orders on the real line and showed that the fractional order entropy
	can better highlight the variations in mathematical, financial and biological series by adjusting the fractional parameter for getting reliable information. \\
	\hspace*{0.2in} This confluence of fractional calculus and entropy appears to be quite challenging as well as helpful in studying the intrinsic details of physical systems. The fractional order entropy conceived by Ubriaco was studied mainly for discrete case but not much work seems to have been done on the continuous case till date as it poses some risk and challenges. So, the focus of our work is to bring to light the limitations of the continuous analogue of Ubriaco entropy (\ref{1.4}), known as the fractional differential entropy (FDE). Henceforth, we worked on addressing the challenges and finding ways to turn it into a desirable entity for studying randomness of continuous distributions associated with complex physical systems. Thus, our current study is sequenced as follows. \\ 
	\hspace*{0.2in} Section 2 is dedicated to examining some statistical properties and bounds of the FDE. Its behavior is then analyzed numerically for different values of $ \alpha $ and compared with the Shannon differential entropy. This section also deals with finding closed form analytical expressions for the entropy formulation of some important and commonly used distributions which may be highly beneficial in reliability and survival analysis works. Lastly, Section 3 applies the fractional entropy in modeling the vertical velocity of wide open channels assuming maximum velocity to act at the surface, ignoring the dip-phenomenon. The estimated velocity is then compared with existing entropy-based models using experimental and field data in order to establish the validity of the model. Throughout the paper, we have taken $ \log(\cdot) $ as the natural logarithm function.
	
	\section{Some properties of FDE}
	For a non-negative absolutely continuous random variable $ X $, the FDE is defined as (Foroghi et al., 2022): 
	\begin{equation}
		H^\alpha (f) = \int_{0}^{\infty} f(x)(-\log f(x))^\alpha~ dx, 
		\label{2.5}
	\end{equation} where $ 0 < \alpha < 1 $. 
	(\ref{2.5}) is the continuous analogue of Ubriaco entropy (\ref{1.4}) which can be seen as the expected value of fractional order negative loglikelihood $ (-\log f(x))^\alpha $, representing the information content of the random variable.
	One can see that for $ \alpha = 1, $ (\ref{2.5}) becomes the differential entropy (\ref{1.1}). \\
	\hspace*{0.2in} Further, we introduce the dynamic versions of the fractional differential entropy. Firstly, we assume that a system started working at time $ t_1 = 0 $ and is found to fail before time $ t_2= t $. Hence, if we consider a random variable $ X_t = [X|X \leq t] $ that denotes the failure time within the interval $ (0,t), $ then the past fractional differential entropy can be defined as:
	\begin{equation*}
		\bar{H}^\alpha(f;t) =  \int_{0}^{t} \frac{f(x)}{F(t)}\Big[-\log \frac{f(x)}{F(t)}\Big]^\alpha~ dx.  
	\end{equation*}
	Here, we can observe that if $ \alpha \to 1, $ then $\bar{H}^\alpha(f;t)$ reduces to the past entropy defined by Di Crescenzo and Longobardi (2002) and when $ t \to \infty, \bar{H}^\alpha(f;t) \to H^\alpha(f)$.
	Similarly, if we know that a system has already survived upto time $ t $, then for the remaining lifetime $ X_t = [X|X\geq t] $ of the used system, the residual fractional differential entropy is defined as:
	\begin{equation*}
		H^\alpha(f;t) = \int_{t}^{\infty} \frac{f(x)}{\bar{F}(t)}\Big[-\log \frac{f(x)}{\bar{F}(t)}\Big]^\alpha~ dx.
	\end{equation*}
	In this case, it is clear that if $ \alpha \to 1, $ then $H^\alpha(f;t)$ reduces to the residual entropy defined by Ebrahimi (1996) and when $ t \to 0, H^\alpha(f;t) \to H^\alpha(f)$. In this current paper, we restrict our study to the fractional differential entropy and leave the dynamic versions for further works.\\
	\hspace*{0.2in} However, direct extension of the fractional entropy proposed by Ubriaco from discrete to continuous case raises some concerns. 
	\begin{itemize}
		\item[(i)] the fractional order differential entropy may assume both positive and negative values for $ f(x) \leq 1 $. 
		\item[(ii)] Another limitation is that for $ f(x) > 1, $ the fractional power of the function $(-\log f(x))$ assumes non-real values. 
	\end{itemize} The situation (ii) is clearly not so desirable. Hence for better applicability of the fractional uncertainty measure, we have restricted our study for absolutely continuous random variables with density function values ranging between $ 0 $ and $ 1.$ This condition can be achieved if we limit the bounds of the random variables or specify the parameters of the distributions. For example, for uniform random variable having pdf $ f(x;A,B) = 1/(B-A) $, we can adjust the upper and lower limits in such a way such that $ B-A \geq 1, ~B > A. $ For exponential variable with pdf $ f(x) = \lambda e^{-\lambda x} $, we can achieve the desired condition if we consider distributions with rate parameter $ \lambda \in (0,1]. $  Similarly, we can achieve the condition $ 0 < f \leq 1 $ for other continuous distributions as well and study the entropy measure for such distributions.\\
	\hspace*{0.2in} In the following, we have evaluated the closed form expressions of fractional order entropy (\ref{2.5}) analytically for ten important families of uni-variate continuous distributions taken from Nadarajah and Zografos (2003), Nanda and Maiti (2007) and Song (2001). The results of which is assumed to be helpful for describing the pdfs commonly useful in many statistical and scientific studies.
	Its essence can be further enhanced by using it to characterize some heavy tailed distributions such as Student's $ t $ distribution with arbitrary degrees of freedom, Cauchy, Cram\'er and Pareto distributions, some of which do not even possess first order moments. 
	Here, we have assumed that the expectation of the fractional power of the negative logarithm of the pdf ($ (-\log f(x))^\alpha ; 0 \leq x \leq C $) exists for the entire support $ [0,C] $ of the variable.
	\begin{enumerate}
		\item \textbf{Uniform distribution}
		$$
		f(x;A,B) = \frac{1}{B-A};~ A \leq x \leq B,~ A \geq 0  $$ 
		$$~and~~	H^\alpha (f) = \big(\log (B-A)\big)^\alpha.$$
		\item \textbf{Normal distribution}
		$$  
		f(x;\mu,\sigma) = \frac{1}{\sqrt{2\pi}\sigma}\exp\Big(-\frac{(x-\mu)^2}{2\sigma^2}\Big); ~x,\mu, \sigma > 0 $$
		$$~and~~ H^\alpha (f) = \frac{1}{2\sqrt{\pi}}\int_{\frac{\mu^2}{2\sigma^2}}^{\infty} e^{-x} x^{-1/2} \big(x - \log(\sqrt{2\pi}\sigma)\big)^\alpha ~dx. $$ 
		\textbf{Particular cases:}
		\begin{equation}
			H^\alpha (f) = \begin{cases}
				\frac{1}{2\sqrt{\pi}}\mathlarger{\mathlarger{\Gamma}}\Big(\alpha + \frac{1}{2},\frac{\mu^2}{2\sigma^2}\Big) & \mbox{;~for } \sigma = \frac{1}{\sqrt{2\pi}}, \mu > 0 \\
				\frac{1}{2\sqrt{\pi}}{\mathlarger{\mathlarger{\Gamma}}}\big(\alpha + \frac{1}{2},1\big) = \frac{1}{2\sqrt{\pi}} E_{\frac{1}{2} - \alpha}(1) = \frac{1}{2\sqrt{\pi}}\varphi_{\alpha - \frac{1}{2}}(1) & \mbox{;~for } \mu^2 = 2\sigma^2 = \frac{1}{\pi} ~,
			\end{cases} \nonumber
		\end{equation} 
		where $$ \Gamma(n,x) = \int_{x}^{\infty} e^{-s}s^{n-1}~ds ~~\&~~ E_m(n) = \int_{1}^{\infty} e^{-nx}x^{-m} dx $$
		are the upper incomplete gamma function and generalized exponential integral, respectively and $ \varphi_{m}(x) = E_{-m}(x) $ is the Misra function (cf. Misra, 1940).
		$ E_m(n) $ can also be expressed as a particular case of the upper incomplete gamma function as follows:
		$$ E_m(n) = n^{m-1}\Gamma (1-m,n). $$  
		\item \textbf{Exponential distribution}
		$$ 
		f(x;\lambda) = \lambda e^{-\lambda x} ;~ x \geq 0 $$
		$$~and~~ H^\alpha (f) = \frac{1}{\lambda} \mathlarger{\Gamma} \big(1+\alpha, \log(1/\lambda)\big). $$
		For $ \lambda = \frac{1}{e}, $ we have 
		$$ H^\alpha (f) = eE_{-\alpha}(1) = e\varphi_{\alpha}(1). $$
		\item \textbf{Pareto II distribution}
		$$
		f(x;k,\sigma) = \frac{k^\sigma \sigma}{(x+k)^{\sigma+1}} ; ~x > 0 ,~ k,\sigma > 0 $$
		$$~and~~ H^\alpha (f) = \Big(\frac{k}{\sigma}\Big)^{\frac{\sigma}{\sigma + 1}} \Big(\frac{\sigma + 1}{\sigma}\Big)^\alpha {\mathlarger{\mathlarger{\Gamma}}}\Big(1+\alpha,\Big(\frac{\sigma}{\sigma + 1}\Big)\log\Big(\frac{k}{\sigma}\Big)\Big).  $$
		\item \textbf{Triangular distribution}
		\begin{equation}
			f(x;\beta) = \begin{cases}
				\frac{2x}{\beta} & \text{;~for } 0 \leq x \leq \beta \\
				\frac{2(1-x)}{\beta} & \text{;~for } \beta < x \leq 1
			\end{cases} \nonumber
		\end{equation} 
		$$~and~~ H^\alpha (f) = \frac{1}{4\cdot 2^\alpha} \mathlarger{\Gamma} \big(1+\alpha,\log(1/4)\big). $$
		\item \textbf{Folded t-distribution}
		$$
		f(x;\gamma) = \frac{2}{\sqrt{\gamma}\beta(\frac{\gamma}{2},\frac{1}{2})} \Big(1 + \frac{x^2}{\gamma}\Big)^{-\frac{\gamma + 1}{2}}; x > 0, \gamma > 0,$$ where $ \beta(m,n) = \int_{0}^{1} u^{m-1} (1-u)^{n-1} du $ is the complete beta function for $ m,n > 0, $ 
		\text{and} \begin{equation*}
			H^\alpha (f) = \frac{2}{(\gamma + 1)\beta(\frac{\gamma}{2},\frac{1}{2})}\int_{0}^{1} (-\log(Ax))^\alpha x^\frac{-1}{\gamma + 1} \bigg(1 - x^\frac{2}{\gamma + 1}\bigg)^{-\frac{1}{2}} dx ,
		\end{equation*} 
		 where $ A = \frac{2}{\sqrt{\gamma}\beta(\frac{\gamma}{2},\frac{1}{2})}.$
		So, for the parameter $ \gamma = 1, A = \frac{2}{\pi} $ and we have: 
		\begin{eqnarray}
			H^\alpha (f) &= & \frac{1}{\beta(\frac{1}{2},\frac{1}{2})}\int_{0}^{1} (-\log(Ax))^\alpha x^\frac{-1}{2} (1 - x)^{-\frac{1}{2}} ~dx \nonumber\\
			&= & \frac{1}{\pi}\int_{0}^{1} (-\log(Ax))^\alpha x^\frac{-1}{2} (1 - x)^{-\frac{1}{2}}~ dx \nonumber\\
			&= & \frac{1}{2} \sum_{k=0}^{\infty} {\frac{1}{2} \choose k} \frac{\mathlarger{\mathlarger{\Gamma}} \big(1+\alpha, \big(\log(\frac{\pi}{2})(k+\frac{1}{2})\big)\big) (k+\frac{1}{2})^{-(1+\alpha)}}{(2/\pi)^{k-\frac{1}{2}}}  \nonumber \\
			&= & \frac{1}{2} \sum_{k=0}^{\infty} {\frac{1}{2} \choose k} \frac{(\log(\frac{\pi}{2}))^{1+\alpha} E_{-\alpha}\big((k+\frac{1}{2})\log(\frac{\pi}{2})\big)}{(2/\pi)^{k-\frac{1}{2}}} \nonumber \\
			&= & \frac{1}{2} \sum_{k=0}^{\infty} {\frac{1}{2} \choose k} \frac{(\log(\frac{\pi}{2}))^{1+\alpha} \mathlarger{\mathlarger{\varphi_{\alpha}}}\big((k+\frac{1}{2})\log(\frac{\pi}{2})\big)}{(2/\pi)^{k-\frac{1}{2}}}, \nonumber
		\end{eqnarray} $~where~~ \mathlarger{{\frac{1}{2} \choose k} = \frac{\prod_{i=0}^{k-1} (\frac{1}{2} - i)}{k!}} ; i \in \mathcal{Z}^+ \cup \{0\}. $ 
		\item \textbf{Cram\'er distribution}
		$$
		f(x;\theta) = \frac{\theta}{2(1+\theta x)^2} ; ~x \geq 0,~ \theta > 0
		$$
		\begin{eqnarray*}
			~and~~	H^\alpha (f) &= & 2^{\alpha - \frac{1}{2}} \Gamma (1+\alpha,\log \sqrt{2})	\\
			&= & 2^{\alpha - \frac{1}{2}} (\log\sqrt{2})^{1+\alpha} E_{-\alpha}(\log\sqrt{2}) \\
			&= & 2^{\alpha - \frac{1}{2}} (\log\sqrt{2})^{1+\alpha} \varphi_{\alpha}(\log\sqrt{2}).
		\end{eqnarray*}
		
		\item \textbf{Cauchy distribution}
		\begin{equation*}
			f(x;\sigma, \mu) = \frac{1}{\pi \sigma(1+(\frac{x-\mu}{\sigma})^2)} ;~ \mu \geq 0,~ \sigma > 0
		\end{equation*}
		\begin{eqnarray*}
			~and~~	H^\alpha (f) &= &  \frac{1}{2}\sqrt{\frac{\sigma}{\pi}}\sum_{k=0}^{\infty} {\frac{1}{2} \choose k} (\pi\sigma)^k\int_{\log A}^{\infty} x^\alpha e^{-(k+\frac{1}{2})x} dx \\
			&= & \frac{1}{2}\sqrt{\frac{\sigma}{\pi}}\sum_{k=0}^{\infty} {\frac{1}{2} \choose k} \frac{(\pi\sigma)^k}{(k+\frac{1}{2})^{1+\alpha}} \mathlarger{\mathlarger{\Gamma}} \Big(1+\alpha,\big(k+\frac{1}{2}\big)\log A\Big) \\
			&= & \frac{1}{2}\sqrt{\frac{\sigma}{\pi}}\sum_{k=0}^{\infty} {\frac{1}{2} \choose k} (\pi\sigma)^k (\log A)^{1+\alpha}\mathlarger{E_{-\alpha}} \Big(\big(k+\frac{1}{2}\big)\log A\Big) \\
			&= & \frac{1}{2}\sqrt{\frac{\sigma}{\pi}}\sum_{k=0}^{\infty} {\frac{1}{2} \choose k} (\pi\sigma)^k (\log A)^{1+\alpha} \mathlarger{\mathlarger{\varphi_{\alpha}}} \Big(\big(k+\frac{1}{2}\big)\log A\Big),
		\end{eqnarray*} where $ A = \pi\sigma\bigg(1+\frac{\mu^2}{\sigma^2}\bigg).$
		\item \textbf{Gamma distribution}
		\begin{equation*}
			f(x;m,n) = \frac{m^n x^{n-1} \exp\{-mx\}}{\Gamma(n)} ;~ m , n > 0 ,~ x > 0
		\end{equation*}
		$$~and~~ H^\alpha (f) = \int_{0}^{\infty} f(x;m,n) [-\log f(x;m,n) ]^\alpha ~dx. $$
		\textbf{Particular case $(n=1)$:}
		$$ H^\alpha (f) = \frac{1}{m} \mathlarger{\Gamma}\big(1+\alpha,\log(1/m)\big). $$
		\item \textbf{Beta distribution}
		\begin{equation*}
			f(x;m,n) = \frac{x^{m-1} (1-x)^{n-1}}{B(m,n)} ;~ m , n > 0 ,~ 0 < x < 1
		\end{equation*}
		\begin{equation*}
			~and~~	H^\alpha (f) =  \begin{cases}
				\mathlarger{\frac{B(m,1)^{\frac{1}{m-1}}(m-1)^\alpha}{m^{1+\alpha}} \mathlarger{\Gamma}\big(1+\alpha,\frac{m}{m-1}\log B(m,1)\big)} & \mbox{;~for } n = 1 \\
				\mathlarger{\frac{B(1,n)^{\frac{1}{n-1}}(n-1)^\alpha}{n^{1+\alpha}} \mathlarger{\Gamma}\big(1+\alpha,\frac{n}{n-1}\log B(1,n)\big)} & \mbox{;~for } m = 1.
			\end{cases}
		\end{equation*}
	\end{enumerate} 
	
	\subsection{Bounds and monotonicity}
	Since it is not possible to derive closed form expressions for the entropy measure for all distributions, we obtain different bounds for the entropy and thereby examine some monotonicity properties of FDE (\ref{2.5}). Before proceeding with the discussion on the bounds, we first need to understand the nature of the entropy function $ h^\alpha(f) = f(x)(-\log f(x))^\alpha $. It is to be noted that the properties are analyzed with respect to $ f $ since the entropy function is considered to be a function of $ f. $ \\ 
	\hspace*{0.2in} In the following proposition, we investigate the point or the value of $ f $ for which the entropy function $ h^\alpha(f) $ assumes the maximum value in the region of our study. The proof is immediate from the case of Ubriaco entropy function ($ s_i(p) = p(i) (-\log p(i))^\alpha $) as discussed in Ubriaco (2009) and hence omitted. 
	\begin{p1}
		The entropy function $ h^\alpha(f) $ attains maximum value at $ f(x) = e^{-\alpha} $, where $ \alpha $ is the order of the entropy function. 
		\label{prop 3.1}
	\end{p1}
\begin{t1}
	\begin{itemize}
		\item[(i)] %
		$ h^\alpha(f) $ decreases and hence the entropy (\ref{2.5}) decreases with increase in $ \alpha $ when $ I(f) = -\log(f(x)) \in (0,1) $, i.e., when $ \frac{1}{e} (= 0.368) < f(x) < 1 $.
		\item[(ii)]Again, $ h^\alpha(f) $ shows an increasing pattern resulting in an increase of entropy with $ \alpha $ when $ I(f) \geq 1 $, i.e., when $ 0 < f(x) \leq 0.368 $. 		
	\end{itemize}
	\label{thm 2.1}
\end{t1}
\begin{figure}[h]
	\centering
	\includegraphics[width=0.5\textwidth]{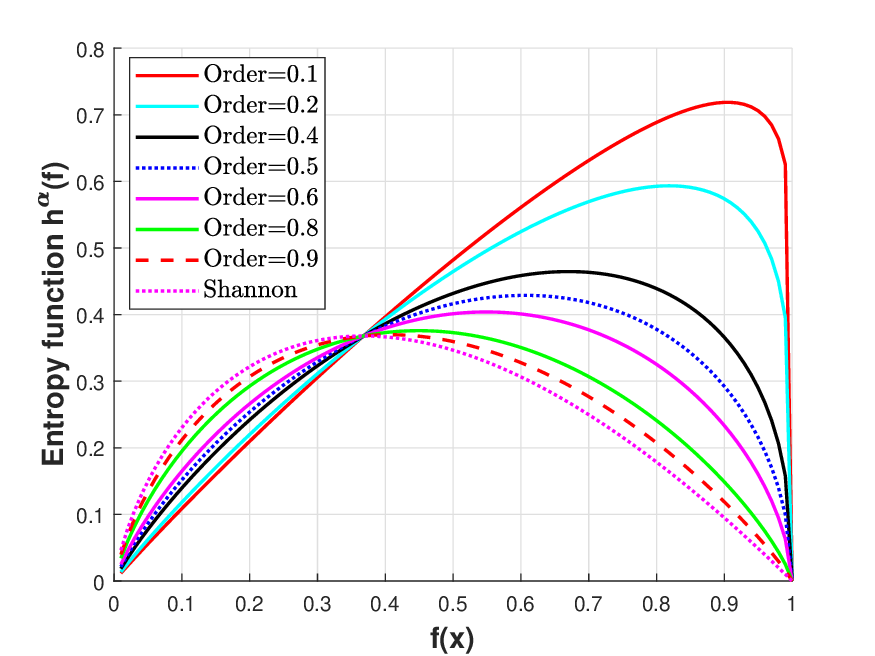} 
	\caption{Fractional order differential entropy function for different orders $ \alpha $}
	\label{ent_order}
\end{figure} 
Fig. \ref{ent_order} portrays the fact that the entropy function $  h^\alpha $ is a non-negative concave function of $ f $ for the region $ 0 < f \leq 1 $. 
Further, we obtain that $ h^\alpha(f) $ increases for $ 0 < f \leq 0.368 $ and decreases when $ 0.368 < f < 1. $ It is interesting to note that when $ I(f) < 0 $ or when $ f(x) > 1 $, then $ h^\alpha(f) $ yields a complex number. So, we ignore such cases since in literature, mostly work related to uncertainty and randomness of a variable is done with real valued measurable functions. \\
\hspace*{0.2in} Further, we have that the differential entropy (\ref{1.1}) is translation invariant but not closed under change of scales or invertible mappings (Shannon, 1948).
In the following proposition, we show that similar property holds for FDE (\ref{2.5}).
\begin{p1}
	Under the linear transformations $ Y = X + c ~and~ Z = mX , m > 0, c \geq 0; $, we have (following Bickel and Lehmann, 2012)
	\begin{eqnarray}
		H^\alpha(f_Y) = H^\alpha(f_X) \label{2.7}\\
		~and~~~~H^\alpha(f_Z) \neq mH^\alpha(f_X).\label{2.8}
	\end{eqnarray}
\end{p1}
\begin{proof}
	The proof for translation invariance (\ref{2.7}) and the relation (\ref{2.8}) follows from the relations $ f_Y(x) = f_X(x-c) ~\text{and}~ f_Z(x) = \frac{1}{m}f_X(\frac{x}{m}), x \in \mathcal{R}.$ 
$\hfill\square$ \\
\vspace{-2em} 
\end{proof}	
In the upcoming propositions, we achieve different bounds of the fractional order entropy in terms of the Shannon differential entropy ($ \ref{1.1} $) when pdf $ f \in (0,1] $. 
\begin{p1}
	For a non-negative random variable~$ X $ with pdf $ f(x) $ 
	\[H^\alpha(f) \leq (H_S(f))^\alpha\] for $ 0 < \alpha \leq 1~ \text{and} ~0 < f \leq 1. $ 
	\label{prop 3.4}
\end{p1}
\begin{proof}
	The result can be achieved from the relation that $ f(x) \leq (f(x))^\alpha $ for $ x > 0,~ $ $0 < f \leq 1 $ and $ 0 < \alpha \leq 1 $ (see Proposition 2 of Xiong et al., 2019).
\end{proof}
\begin{p1}
	If ~$ X $ has a bounded support $ [0,b] $ then for $ 0 < \alpha \leq 1~ \text{and} ~0 < f \leq 1. $ 
	\begin{equation}
		H^\alpha(f) \geq A(\alpha)e^{H_S(f)} \label{2.9}
	\end{equation} where $ A(\alpha) = \exp\left[\int_{0}^{b} f(x)\log(f(x)[-\log f(x)]^\alpha)\,dx\right] $.
\end{p1}
\begin{proof}
	Using the log-sum inequality, we have
	\begin{equation}
		\int_{0}^{\infty} f(x)\log \frac{f(x)}{f(x)[-\log f(x)]^\alpha}dx \geq \log \frac{1}{\int_{0}^{\infty}f(x)[-\log f(x)]^\alpha dx} = \log \frac{1}{H^\alpha (f)}. \label{2.10}
	\end{equation}
	Moreover, the L.H.S of (\ref{2.10}) can be written as
	\begin{equation}
		\int_{0}^{\infty} f(x)\log \frac{f(x)}{f(x)[-\log f(x)]^\alpha}dx = -H_S(f)-\int_{0}^{\infty}f(x)\log (f(x)[-\log f(x)]^\alpha)dx.
		\label{2.11}
	\end{equation}
	Now, using the result of Proposition \ref{prop 3.1}, we get that \begin{equation}
		\log (f(x)[-\log f(x)]^\alpha) \leq \alpha\log \frac{\alpha}{e}. \label{2.12}	
	\end{equation} 
	Hence, if $ X $ has a bounded support $ [0,b] $, then from (\ref{2.10}) and (\ref{2.12}), we get that \begin{equation} A(\alpha) = \exp[\int_{0}^{b} f(x)\log(f(x)[-\log f(x)]^\alpha)\,dx] \leq b \log (e^{-\alpha}\cdot\alpha^\alpha), \label{2.13} 
	\end{equation}
	where we can find that the function $ \log(e^{-\alpha}\cdot\alpha^\alpha) $ is monotonically decreasing w.r.t $ \alpha \in (0,1] $ and hence $ A(\alpha) < \infty $ and bounded. Therefore, from (\ref{2.11}) and (\ref{2.13}), we get (\ref{2.9}).
\end{proof}
\begin{p1}
	Let $ X $ be bounded with support $ [0,b]~,b > 0 $. Then for $ 0 < \alpha \leq 1~ \text{and} ~0 < f \leq 1, $ 
	\begin{equation}
		H^\alpha(f) \leq b^{1-\alpha}(H_S(f))^\alpha.
	\end{equation}
\end{p1}
\begin{proof}
	The proof follows from Proposition \ref{prop 3.4} and Jensen's inequality for integrals (see Proposition 3.2 of Di Crescenzo et al., 2021).
	$\hfill\square$ \\
	\vspace{-1em} 
\end{proof}
In the next two propositions, we obtain some more bounds for the FDE with respect to the negative loglikelihood function and some integrals over a bounded domain $ [0,b]. $
\begin{p1}
	Let $ X $ be bounded in $ L^1 $ space having support $ [0,b], ~b > 0, $ then 
	\begin{equation}
		H^\alpha(f) \leq \int_{0}^{b} I^\alpha(f) ~dx
	\end{equation} for $ 0 < f \leq 1, $ where $ I^\alpha(f) = [-\log f(x)]^\alpha $ is the information gain function of FDE.
\end{p1}
\begin{proof}
	This result follows from Holder's inequality which states that for measurable real valued functions $ f_1 $ and $ f_2 $ defined on a domain $ D $ bounded in $ L^1 $ space with closed support $ [0,b] $, we have that
	\[\left|\int f_1(x)f_2(x) d\mu(x)\right| \leq \left(\int |f_1(x)| d\mu(x)\right)\left(\int |f_2(x)| d\mu(x)\right).\]
	Here, we have taken $ f_1(x) = f(x) $, $ f_2(x) = (-\log f(x))^\alpha .$
\end{proof}

\begin{p1}
	For a random variable $ X $ with bounded support $ [0,b] $ and $ H^\alpha(f) < \infty $, we have
	\begin{itemize}
		\item[(i)] $
		\int_{0}^{b} f(x)\big(1 -f(x)\big)^\alpha\,dx \leq H^\alpha(f) \leq \int_{0}^{b} f(x)\bigg(\frac{1}{f(x)}-1\bigg)^\alpha\,dx ;
		$
		\item[(ii)] $ H^\alpha(f) \leq b\big(\frac{\alpha}{e}\big)^\alpha $
	\end{itemize}
	if $ 0 < f \leq 1~ \text{and} ~ \alpha > 0. $ 
\end{p1}
\begin{proof}
	\begin{itemize}
		\item[(i)] The inequalities in $(i)$ can be achieved by using the relation 
		\[1-\frac{1}{x} \leq \log x \leq x - 1\] for $ x > 0. $ 
		\item[(ii)] For the second part, we have used the fact that $ f(x)[-\log f(x)]^\alpha $ is non-negative and concave for $ 0 < f(x) \leq 1 ~\forall~ \alpha \in (0,1] $ such that $$ f(x)[-\log f(x)]^\alpha \leq [-\log a]^{\alpha - 1}[\alpha a - f(x)(\alpha + \log a)], $$ with $ \alpha \in (0,1] ~\text{and}~ a \in (0,1] $ (see Proposition 3.3 of Di Crescenzo et al., 2021). \\
		Therefore, from equation (\ref{2.5}), we have
		$$ H^\alpha(f) \leq [-\log a]^{\alpha - 1}[\alpha a b - (\alpha + \log a)]. $$
		By substituting $ a = e^{-\alpha}, $ we get the final expression given by $(ii)$.	
		$\hfill\square$ \\
	\end{itemize}	
\vspace{-3em}  
\end{proof}
In the following proposition, we inspect the non-additivity property of FDE and used it to obtain a bound for the entropy measure for bivariate distribution.
\begin{p1}
	For two non-negative absolutely continuous independent random variables $ X $ and $ Y $,
	\begin{equation}
		H^\alpha(f_X,f_Y) \leq H^\alpha(f_X) + H^\alpha(f_Y),~ 0 < \alpha \leq 1
	\end{equation} for $ 0 < f \leq 1 $ with equality only if $ \alpha = 1. $
\end{p1}
\begin{proof}
	The proof follows from the relation that
	\[\iint f(x,y)\,dx\,dy = \int f_X(x)\,dx \int f_Y(y)\,dy\] for independent random variables~$X$ and ~$ Y $ and then using the inequality $$ (x+y)^\alpha \leq x^\alpha + y^\alpha;~~0 < \alpha \leq 1,~~ x,y > 0, $$ with equality only if $ \alpha = 1. $ $\hfill\square$ \\
	\vspace{-2em} 
\end{proof}
In reliability modeling works, the significance of a entropy measure is decided by its values. Higher the value of the entropy, more information is associated with the distribution, making it more suitable to describe a particular system or a random variable. Driven by this property of entropy function, we investigate the effectiveness of FDE $ H^\alpha(f) $ over Shannon entropy $ H_S(f) $. It is evident from Theorem \ref{thm 2.1} that $ H^\alpha(f) $ contains more information for some $ \alpha $ close to zero whereas $ H_S(f) $ gives the minimum entropy for those distributions whose $ f(x) \in (0.368,1). $

\begin{table}[h]
	\small
	\centering
	\caption{Some important continuous distributions}
	\begin{tabular}{| l | l |}
		\hline
		\textbf{Distributions} & ~~~~ \textbf{Probability Density Function} \\
		\hline \hline
		Weibull & ~~~~$f(x;a,b) =  \frac{b}{a}(\frac{x}{a})^{b-1}\exp(-(\frac{x}{a})^b),~ x \geq 0 $ \\
		Uniform & ~~~~$ f(x;A,B) = \frac{1}{B-A},~ A \leq x \leq B $ \\
		\hline \hline
		Standard Normal & ~~~~ $f(x) = \frac{1}{\sqrt{2\pi}}\exp(-\frac{x^2}{2}),~ 0 < x < \infty $ \\
		Generalized Pareto & ~~~~$ f(x;k,\sigma,\theta) = \big(\frac{1}{\sigma}\big) \bigg(1 + \frac{k(x-\theta)}{\sigma}\bigg)^{-1-\frac{1}{k}}$,~ $\theta < x $ for $ k > 0 $ \\
		Exponential & ~~~~$ f(x;\lambda) = \lambda e^{-\lambda x},~ x \geq 0, ~\lambda > 0 $ \\
		Finite Range & ~~~~$ f(x;a,\theta) = \frac{a}{x} \big(\frac{x}{\theta}\big)^a, ~0 \leq x \leq \theta,~ a,\theta > 0. $ \\
		\hline
	\end{tabular}	
	\label{tab:dist2}	
\end{table}	

\begin{table}[th]
	\small
	\caption{Fractional order entropy values for continuous distributions}	
	\begin{tabular}{|l | c  c  c  c  c  c |}
		\hline
		\textbf{Distributions} &~ $ \alpha = 1 $ &~ $ \alpha = 0.9 $ &~ $ \alpha = 0.8 $ &~ $ \alpha = 0.6 $ &~ $ \alpha = 0.4 $ &~ $ \alpha = 0.2 $ \\
		\hline \hline
		Weibull ($ a=1, b=2 $) &~ 0.2300 &~ 0.2389 &~ 0.2499 &~ 0.2790 &~ 0.3201 &~ 0.3768 \\
		\hline
		Uniform ($ A=0,B=2 $) &~ 0.6931 &~ 0.7190 &~ 0.7459 &~ 0.8026 &~ 0.8636 &~ 0.9293 \\
		\hline \hline
		Standard Normal &~ 1.4183 &~ 1.3580 &~ 1.3030 &~ 1.2068 &~ 1.2161 &~ 1.0579 \\
		\hline
		GPD ($ k = 1, \sigma = 2, \theta = 2 $ ) &~ 1.4025 &~ 1.3229 &~ 1.2485 &~ 1.1136 &~ 0.9953 &~ 0.8914 \\
		\hline
		Exponential($ \lambda = 1 $) &~ 1.0000 &~ 0.9618 &~ 0.9314 &~ 0.8935 &~ 0.8873 & 0.9182 \\
		\hline
		Finite Range ($ p = 2,\theta = 10 $) &~ 1.4071 &~ 1.3520 &~ 1.2912 &~ 1.2074 &~ 1.1257 &~ 1.0537 \\
		\hline	
	\end{tabular}
	\label{tab:entvalue}
\end{table}	

\subsection{Numerical comparison}
According to the principles of information theory, the entropy measure that yields a higher value for a specific distribution, when fitted to a given dataset of a physical system, is considered more effective in representing the system's uncertainty or information content.
In this subsection, in an attempt to resolve the quest for finding those distributions which would go well with Theorem \ref{thm 2.1} $(i)$, where $ H^\alpha(f) $ can be considered more efficient over $ H_S(f) $, giving maximum entropy as $ \alpha $ approaches zero whereas attaining minimum value at $ \alpha = 1 $, we have considered some important continuous distributions popularly found in reliability and survival analysis studies enlisted in Table \ref{tab:dist2}.\\ 
\hspace*{0.2in} Referring to the numerically computed entropy values for different values of $ \alpha $ for those particular distributions provided by Table \ref{tab:entvalue}, one can find that for Weibull and Uniform distributions under certain parameters, minimum entropy value is obtained at $ \alpha=1$ and maximum is attained for $ \alpha $ close to zero. However, an opposite trend is seen for the remaining distributions (Standard Normal, GPD, Exponential and Finite Range), where we get the lowest value for $ \alpha $ near zero and the value increases with increase in $ \alpha $. This changing behavior of $ H^\alpha(f) $ with different distributions are a consequence of the property of $ h^\alpha(f) $ mentioned in Theorem \ref{thm 2.1}, thereby imparting more flexibility to the application of FDE than the Shannon differential entropy in checking fitness of a model describing a system.  
\section{Modeling of one-dimensional velocity distribution through FDE}
Complex systems such as open-channel hydraulics display a certain extent of chaotic behaviour. Therefore, the measurement or prediction of the distribution of the hydraulic variables such as velocity, suspended sediment concentration and sheer stress is quite significant for flood control measures and water quality management. For describing the one-dimensional flow velocity distribution through open channels, the linear law was developed to represent near bed flow velocity (viscous sublayer), Prandtl-von Karman universal velocity distribution law for middle region of the channels and power law for near water surface velocity (see Sarma et al., 1983 and references therin). But, these laws failed to give one single equation for the flow velocity along the entire vertical region from channel bed to the surface of the channel. 
Further, the modeling of velocity distribution in unsteady turbulent flows through conventional deterministic models poses quite a lot of challenges such as computational complexity and limited availability of data or incomplete information due to human errors and sampling difficulty. To overcome these limitations, probabilistic concepts were introduced in hydraulic studies to derive entropy based velocity distribution expressions for open channels as entropy neither depend on the symmetry of the distribution nor require metric data for its computation. Moreover, it gives accurate results even with small samples. This in turn reduces the computational time and cost. These advantages of the entropy-based methods for estimating the distribution of the flow velocity of open channels attracted a lot of researchers to develop their models based on different kinds of entropy to present better and simpler estimation methods. The existing one-dimensional vertical velocity models, assuming maximum surface velocity, are based on Shannon entropy (Chiu, 1987), Tsallis entropy (Singh and Luo, 2011), R$\acute{e}$nyi entropy (Kumbhakar and Ghoshal, 2016) and fractional Machado entropy (Kumbhakar and Tsai, 2023). \\
\hspace*{0.2in} So far, to the best of our knowledge, the FDE due to Ubriaco has not been taken into consideration to derive the vertical distribution of flow velocity in open-channels. Hence, it might be appealing to explore the use of the FDE to frame an entropy theory for determining the one-dimensional velocity distribution of turbulent flows in wide rectangular open channels, such as rivers, along the vertical depth by using the principle of maximum entropy (POME) proposed by Jaynes (Jaynes, 1957). The POME states that the most preferred or the least-biased probability distribution is the one which maximizes the entropy, with respect to some constraints based on some prior information associated with the random variable.  In practical situations, due to the constraints, the probability distribution that maximizes the entropy may not be uniform. Hence, having a maximum entropy implies getting a distribution close to the uniform distribution while satisfying the constraint equations. 
Therefore, by using the method of calculus of variations, we maximize our entropy
\begin{equation}
	H^\alpha(f) = \int_{0}^{1} f(\hat{\nu})[-\log f(\hat{\nu})]^\alpha~d\hat{\nu}, \label{apl_1} 
\end{equation} where the random variable $ \mathcal{\hat{V}} $ with pdf $ f(\hat{\nu}) $ attains the dimensionless velocity $ \hat{\nu} = \frac{\nu}{\nu_{max}}. $ It normalizes the vertical velocity $ \nu $ of the stream with $ \nu_{max} $ being the maximum velocity at a particular height from the channel bed, subject to a specific set of constraints given by
\begin{flalign*}
	&(1)~~ \int_{0}^{1} f(\hat{\nu})~d\hat{\nu} = 1 & 
	&(3)~~ \int_{0}^{1} \hat{\nu}^2f(\hat{\nu})~d\hat{\nu} = \delta\hat{\nu}_m^2 & \\
	&(2)~~ \int_{0}^{1} \hat{\nu} f(\hat{\nu})~d\hat{\nu} = \hat{\nu}_m & 
	&(4)~~ \int_{0}^{1} \hat{\nu}^3f(\hat{\nu})~d\hat{\nu} = \beta\hat{\nu}_m^3~. &	
\end{flalign*}
\hspace*{0.2in} Here, $ \hat{\nu}_m $ is the vertical cross-sectional mean flow velocity, $ \delta $ and $ \beta $ are the coefficients of momentum and energy distribution, respectively. 
Constraint $(1)$ is called the total probability law or the complete normalization condition.
Constraints (2), (3) and (4) are the $ 1^{st}, 2^{nd} $ and $ 3^{rd} $ order moments of the random variable $ \mathcal{\hat{V}} $, respectively. These moments represent the conservation of mass, momentum and energy, respectively, in the process of flow through the cross-section of a particular channel. Since it is assumed that the effects of momentum and energy on the flow velocity is negligible and the law of conservation of mass is obeyed by every flow, hence we consider only the first two constraints to construct the Lagrangian in order to estimate the distribution with maximum entropy.
The Lagrangian function $ L $ can be written as:
\begin{equation}
	L(f,\hat{\nu}) = \int_{0}^{1} f(\hat{\nu})[-\log f(\hat{\nu})]^\alpha~d\hat{\nu} -  a\Big\{\int_{0}^{1} f(\hat{\nu})~d\hat{\nu} - 1\Big\} - b\Big\{\int_{0}^{1} \hat{\nu} f(\hat{\nu})~d\hat{\nu} - \hat{\nu}_m\Big\}, \label{apl_2}
\end{equation}
where $ a $ and $ b $ are the Lagrange multipliers.
We can obtain the pdf $ f(\hat{\nu}) $ of $ \mathcal{\hat{V}} $ by ignoring the integration sign on the right hand side of (\ref{apl_2}) and solving the Euler-Lagrange equation considering $ f(\hat{\nu}) $ and $ \hat{\nu} $ as the dependent and independent variables, respectively, given by
\begin{equation*}
	\frac{\partial L}{\partial f(\hat{\nu})} - \frac{d}{d\hat{\nu}}\Big(\frac{\partial L}{\partial f(\hat{\nu})'}\Big) = 0,
\end{equation*} we have
\begin{equation}
	\frac{\partial f(\hat{\nu})[-\log f(\hat{\nu})]^\alpha}{\partial f(\hat{\nu})} ~-~ a~ \frac{\partial f(\hat{\nu})}{\partial f(\hat{\nu})} ~-~ b~\frac{\partial \hat{\nu} f(\hat{\nu})}{\partial f(\hat{\nu})} = 0.\label{apl_3} 
\end{equation}
Simplifying (\ref{apl_3}) gives us
\begin{equation}
	(-\log f(\hat{\nu}))^\alpha - \alpha(-\log f(\hat{\nu}))^{\alpha - 1} - a - b\hat{\nu} = 0. \label{apl_4}
\end{equation}
 If we consider $ a + b\hat{\nu} = x(\hat{\nu}) $ and $ -\log f(\hat{\nu}) = y(\hat{\nu}) $, then (\ref{apl_4}) can be expressed as:
\begin{equation}
	y(\hat{\nu})^{\alpha} - \alpha y(\hat{\nu})^{\alpha - 1} - x(\hat{\nu}) = 0. \label{apl_5} 
\end{equation}
\hspace*{0.2in} We find that eqn. (\ref{apl_5}) is an implicit non-linear equation in $ y(\hat{\nu}) $ which can be solved numerically.
To obtain an explicit expression for $ y(\hat{\nu}) $ in terms of $ x(\hat{\nu}), $ we have taken a particular value of $ \alpha $ (say, $ \frac{1}{2}). $ 
 Therefore, for $ \alpha = \frac{1}{2} $, we get a quadratic equation in $ y(\hat{\nu}) $ given by
\begin{equation}
	y(\hat{\nu})^2 - y(\hat{\nu})(1+x(\hat{\nu})^2) + \frac{1}{4} = 0. \label{apl_6} 
\end{equation}
Solving (\ref{apl_6}), we get
\begin{equation}
	y(\hat{\nu}) = \frac{(1+x(\hat{\nu})^2)\pm x(\hat{\nu}) \sqrt{x(\hat{\nu})^2 + 2}}{2} . \nonumber
\end{equation}
 This implies that
\begin{eqnarray}
	f(\hat{\nu}) &=& \exp\Bigg[-\frac{(1+(a+b\hat{\nu})^2) + (a+b\hat{\nu}) \sqrt{(a+b\hat{\nu})^2 + 2}}{2}\Bigg] \label{apl_7} \\ 
	\text{or} ~~f(\hat{\nu}) &=& \exp\Bigg[-\frac{(1+(a+b\hat{\nu})^2) - (a+b\hat{\nu}) \sqrt{(a+b\hat{\nu})^2 + 2}}{2}\Bigg]. \label{apl_8} 
\end{eqnarray}
\hspace*{0.2in} Substituting the value of $ f(\hat{\nu}) $ from (\ref{apl_7}) in (\ref{apl_1}), we obtain the maximum entropy of $ f(\hat{\nu}) $ of $ \mathcal{\hat{V}} $.
Therefore, we have the FDE which can be computed as follows:
\begin{align}
	H^\alpha(f) = \int_{0}^{1} \exp\Bigg[-\frac{(1+(a+b\hat{\nu})^2) + (a+b\hat{\nu}) \sqrt{(a+b\hat{\nu})^2 + 2}}{2}\Bigg]\times \nonumber\\
	\Bigg(\frac{(1+(a+b\hat{\nu})^2) + (a+b\hat{\nu}) \sqrt{(a+b\hat{\nu})^2 + 2}}{2}\Bigg) ~d\hat{\nu}. \label{apl_9} 
\end{align}
\subsection{Estimation of the velocity distribution}
\begin{figure}
	\centering
	\includegraphics[width=0.45\textwidth]{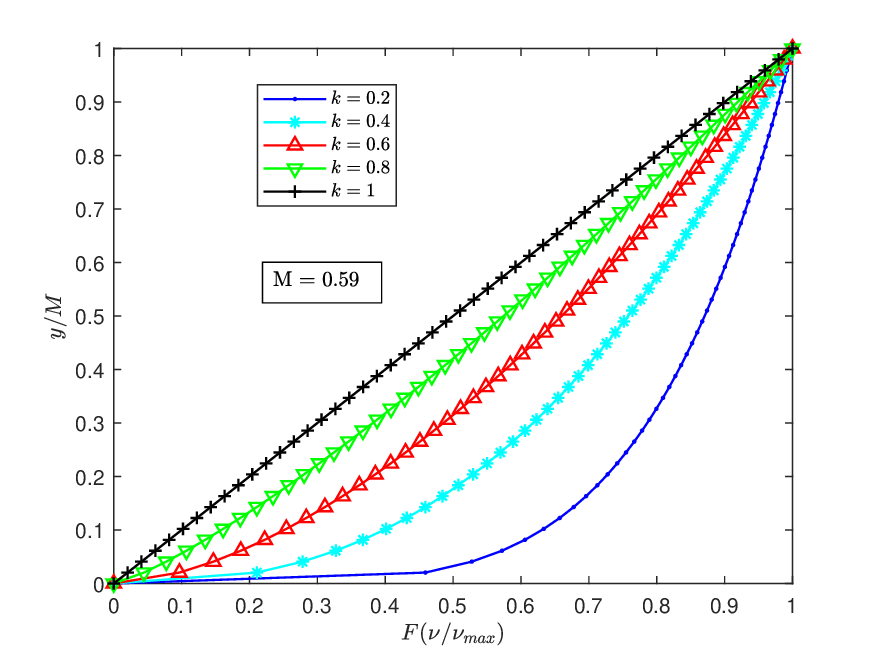}
	\caption{Variation of cdf with respect to the fitting parameter $ k $}
	\label{cdf_k}
\end{figure}

For deriving a general expression for one-dimensional velocity distribution along the vertical column of an open channel, firstly we need to find an appropriate relation between the cdf of normalized dimensionless flow velocity $ \mathcal{\hat{V}} $ and the normalized height $ y/M $. Here $ M $ symbolizes the total height of the water channel and $ y $ is the vertical distance from channel bed. 
This spatial relation of cdf $ F(\hat{\nu}) $ can be hypothesized by taking into account that the vertical flow velocity increases monotonically from 0 to $ \hat{\nu}_{max} $ as height $ y $ increases from $ y=0 $ at the channel bed to a maximum level $ y=M $ at the water surface and, also considering the geometry of the channel and the characteristics of suspended sediment particles throughout the flow domain which affect the gradient of velocity with respect to the vertical height. Further, the cdf should be continuous and differentiable within the vertical height range of the channel $(0,M)$ and its values should lie between 0 and 1. The density function can then be derived by differentiating $ F(\hat{\nu}) $ with respect to $ \hat{\nu}. $ Thus, using the chain rule of differentiation, we obtain  
\begin{equation}
	f(\hat{\nu}) = \frac{dF(\hat{\nu})}{d\hat{\nu}} = \frac{dF(\hat{\nu})}{dy}\cdot \frac{dy}{d\hat{\nu}} .\label{est1} 
\end{equation}
Since, $ f(\hat{\nu}) \geq 0 $ and $ \frac{dy}{d\hat{\nu}} \geq 0, $ from eqn. (\ref{est1}) one can see that $ \frac{dF(\hat{\nu})}{dy} \geq 0. $ In other words, the cdf will be an increasing function of the space domain $ y/M. $ \\ 
\hspace*{0.2in} Considering the above stated conditions, the spatial cdf relation can be defined as:
\begin{equation}
	F(\hat{\nu}) = \left(\frac{y}{M}\right)^k,
	\label{est2}
\end{equation} where $ k \in [0,1] $ is taken as the fitting parameter representing the sediment particle characteristics, which can be determined by least square method using a given data set (Cui and Singh, 2012). Smaller $ k $ values will represent less impact of the suspended particles on the velocity which will lead to less rate of increase of cdf with height. Larger $ k $ values will represent greater impact, leading to rise in rate of increase in cdf with height. This influence of $ k $ values on the hypothetical cdf (\ref{est2}) is illustrated in Fig. \ref{cdf_k}. We can observe that $ k = 1 $ corresponds to a linear relation of the cdf $ F(\hat{\nu}) $ with the normalized height $ y/M. $ This condition is usually prevalent in clear water channels or in steady flows, where the velocity $ \nu $ increases uniformly with height. \\
\hspace*{0.2in} Now, the FDE based cdf relation with respect to the velocity domain can be determined from eqn. (\ref{apl_7}) or (\ref{apl_8}). Since, we know that the cdf is an increasing function of $ \hat{\nu}, $ we choose eqn. (\ref{apl_7}) as the least biased pdf of $ \mathcal{\hat{V}}. $
Therefore, using eqn. (\ref{apl_7}), the cdf in terms of $ \hat{\nu} $ can be obtained as:
\begin{eqnarray}
	F(\hat{\nu}) &=& \int_{0}^{\hat{\nu}} f(u) du \nonumber\\ 
	&=& \int_{0}^{\hat{\nu}} \exp\Bigg[-\frac{(1+(a+bu)^2) + (a+bu) \sqrt{(a+bu)^2 + 2}}{2}\Bigg] du  
	\label{est3} 
\end{eqnarray}
\hspace*{0.2in} 
Using the substitution $ (a+bu)^2 = t $ in eqn. (\ref{est3}) with the assumption $ |t| < 1, $ we get
\begin{equation}
	F(\hat{\nu}) = \frac{e^{-1/2}}{2b}\int_{a^2}^{(a+b\hat{\nu})^2}t^{-1/2}\exp\Big[-\frac{t+\sqrt{2t+t^2}}{2}\Big] dt \label{est4}
\end{equation}
\hspace*{0.2in}
 Now, using the power series expansion of $ e^{-x} = \sum_{i=0}^{\infty} (-1)^i \frac{x^i}{i!}, ~ \forall x \in \mathcal{R} $ in (\ref{est4}), we get
\begin{eqnarray*}
	F(\hat{\nu}) &=& \frac{e^{-1/2}}{2b}\int_{a^2}^{(a+b\hat{\nu})^2} t^{-1/2}\sum_{i=0}^{\infty}\frac{(-1)^i}{i!}\Bigg(\frac{t+\sqrt{2t+t^2}}{2}\Bigg)^i dt.	
\end{eqnarray*}
\hspace*{0.2in} Again, using binomial expansion formulas given by $ (x+y)^i = \sum_{j=0}^{i}{i \choose j} x^{i-j}y^j,~i\in \mathcal{N} $ and $ (1+x)^i = \sum_{j=0}^{\infty}{i \choose j} x^j,~i\in \mathcal{R} $ for $ |x| < 1 $, we obtain
\begin{eqnarray}
	F(\hat{\nu}) &=& \frac{e^{-1/2}}{2b}\sum_{i=0}^{\infty}\frac{(-1)^i}{2^i i!}\int_{a^2}^{(a+b\hat{\nu})^2} t^{-1/2}\sum_{j=0}^{i}{i \choose j}t^{i-j}(2t+t^2)^{j/2} dt	\nonumber\\
	&=& \frac{e^{-1/2}}{2b}\sum_{i=0}^{\infty}\frac{(-1)^i}{2^i i!}\sum_{j=0}^{i}{i \choose j}\int_{a^2}^{(a+b\hat{\nu})^2} t^{-1/2}t^{i-j}(2t)^{j/2}\Big(1+\frac{t}{2}\Big)^{j/2} dt; ~\Big|\frac{t}{2}\Big| < 1 \nonumber\\ 
	&=& \frac{e^{-1/2}}{2b}\sum_{i=0}^{\infty}\sum_{j=0}^{i}\frac{(-1)^i}{2^i i!}{i \choose j}2^{j/2} \int_{a^2}^{(a+b\hat{\nu})^2} t^{-\frac{1}{2}+i-\frac{j}{2}} \sum_{k=0}^{\infty}{j/2 \choose k}\Big(\frac{t}{2}\Big)^k ~dt \nonumber\\
	&=&
	\frac{e^{-1/2}}{b}\sum_{i=0}^{\infty}\sum_{j=0}^{i}\sum_{k=0}^{\infty}{i \choose j} {j/2 \choose k}\frac{(-1)^i ~2^{j/2 - i - k}}{i!}\cdot\frac{(a+b\hat{\nu})^{1+2(i+k)-j}-a^{1+2(i+k)-j}}{1+2(i+k)-j}.\nonumber \\
	&& \label{est5}
\end{eqnarray}
\hspace*{0.2in} For computational convenience, we assume $ |a| < 1 ~\text{and}~ |a+b| < 1 $ so that we can approximate (\ref{est5}) up to two terms and ignore the higher power terms of $ (a+b) $ and $ a. $ 
Through this, we get an explicit expression of $ F(\hat{\nu}) $ for the case $(i=0,1; k=0) $, expressed as:
\begin{equation}
	F(\hat{\nu}) = -\frac{e^{-1/2}}{2\sqrt{2}b}\left[(a+b\hat{\nu})^2-a^2\right]. \label{est6}
\end{equation}
 \hspace*{0.2in} Therefore, by comparing the proposed spatial cdf (\ref{est2}) with the FDE based cdf (\ref{est6}), we get the closed form expression for the normalized velocity distribution given as:
\begin{equation}
	\hat{\nu} = \frac{1}{b}\Bigg[-a\pm\sqrt{a^2 - 2\sqrt{2}be^{1/2}\left(\frac{y}{M}\right)^k} \Bigg] 
	\label{est7} 
\end{equation}
\hspace*{0.2in} The expression (\ref{est7}) represents the FDE-based dimensionless normalized velocity distribution for the centre region of wide open channels where maximum velocity is assumed to occur at the channel surface.

\subsection{Determination of Lagrange multipliers}
The FDE-based velocity distribution expression defined by eqn. (\ref{est7}) contains three unknown terms $ a $, $ b $ and $ k, $ which need to be determined. The Lagrange multipliers $ a $ and $ b $ can be determined by substituting $ f(\hat{\nu}) $ obtained from eq. (\ref{apl_7}) in the constraints (1) and (2). 
By putting the value of $ f(\hat{\nu}) $ in the constraint (1) and by using the transformation $ (a+b\hat{\nu})^2 = t $, we have
\begin{eqnarray}
	L.H.S &=& \int_{0}^{1} f(\hat{\nu})~d\hat{\nu} \nonumber\\
	&=& \frac{e^{-1/2}}{2b}\int_{a^2}^{(a+b)^2}t^{-1/2}\exp\Big[-\frac{t+\sqrt{2t+t^2}}{2}\Big] dt. \label{lag1} 
\end{eqnarray}
\hspace*{0.2in} Following similar steps used in the determination of the FDE-based cdf expression represented by eqn. (\ref{est6}), we come up with the final expression given by:
\begin{equation*}
	L.H.S =  \frac{e^{-1/2}}{b}\sum_{i=0}^{\infty}\sum_{j=0}^{i}\sum_{k=0}^{\infty}{i \choose j}{j/2 \choose k}\frac{(-1)^i 2^{\frac{j}{2}-k-i}}{i!} \Bigg[\frac{(a+b)^{1+2i-j+2k}-a^{1+2i-j+2k}}{1+2i-j+2k}\Bigg]. 
\end{equation*}
\hspace*{0.2in} Now, L.H.S = 1 yields:
\begin{equation}
	\sum_{i=0}^{\infty}\sum_{j=0}^{i}\sum_{k=0}^{\infty}{i \choose j}{j/2 \choose k}\frac{(-1)^i 2^{\frac{j}{2}-k-i}}{i!} \Bigg[\frac{(a+b)^{1+2i-j+2k}-a^{1+2i-j+2k}}{1+2i-j+2k}\Bigg] = be^{1/2}. \label{lag2} 
\end{equation}
\hspace*{0.2in} By assuming $ |a| < 1 ~\text{and}~ |a+b| < 1 $ to approximate (\ref{lag2}) up to two terms and ignoring the higher power terms of $ (a+b) $ and $ a, $ we get a linear equation in $ a $ and $ b $. Hence, for the case $(i=0,1; k=0) $, we get
\begin{equation}
	b+2a = -2\sqrt{2}e^{1/2}. \label{lag3}
\end{equation}
\hspace*{0.2in}Similarly, by putting the value of $ f(\hat{\nu}) $ in the constraint (2) and on using the transformation $ a+b\hat{\nu} = t $, we have
\begin{eqnarray}
	L.H.S &=& \int_{0}^{1} \hat{\nu}f(\hat{\nu})~d\hat{\nu} \nonumber\\
	&=& \frac{1}{b}\int_{a}^{a+b} \Big(\frac{t-a}{b}\Big)
	\exp\Big[-\Big(\frac{1+t^2}{2}\Big)\Big]\cdot\exp\Big[-\frac{t\sqrt{2+t^2}}{2}\Big]~dt \nonumber \\
	&=& \frac{e^{-1/2}}{b^2}\Bigg[\int_{a}^{a+b}t\cdot\exp\Big[-\frac{t^2+\sqrt{2t^2+t^4}}{2}\Big]~dt - a\int_{a}^{a+b}\exp\Big[-\frac{t^2+\sqrt{2t^2+t^4}}{2}\Big]~dt\Bigg]. \nonumber\\
	\text{Now},~ I_1 &=& \int_{a}^{a+b}t\cdot\exp\Big[-\frac{t^2+\sqrt{2t^2+t^4}}{2}\Big]~dt \nonumber\\
	&=& \sum_{i=0}^{\infty}\frac{(-1)^i}{2^i i!}\sum_{j=0}^{i} {i \choose j}\int_{a}^{a+b}t^{1+2(i-j)}(2t^2+t^4)^{j/2}~dt. \nonumber 
\end{eqnarray}
\hspace*{0.2in} Let us assume $ |a|<1 $ and $|a+b|<1 \implies |t|< 1 $ since $ 0 \leq \hat{\nu} \leq 1. $
\begin{eqnarray}
	\implies I_1 &=& \sum_{i=0}^{\infty}\frac{(-1)^i}{2^i i!}\sum_{j=0}^{i} {i \choose j}\int_{a}^{a+b}t^{1+2(i-j)+j}\cdot2^{j/2}\Big(1+\frac{t^2}{2}\Big)^{j/2}~dt ;~ \Big|\frac{t^2}{2}\Big| < 1 \nonumber \\
	&=& \sum_{i=0}^{\infty}\frac{(-1)^i}{2^i i!}\sum_{j=0}^{i} {i \choose j}\sum_{k=0}^{\infty} {j/2 \choose k} 2^{j/2 - k}\int_{a}^{a+b}t^{1+2i-j+2k}~dt \nonumber \\
	&=& \sum_{i=0}^{\infty}\sum_{j=0}^{i}\sum_{k=0}^{\infty}{i \choose j} {j/2 \choose k}\frac{(-1)^i ~2^{j/2 - i - k}}{i!}\Bigg[\frac{(a+b)^{2(1+i+k)-j} - a^{2(1+i+k)-j}}{{2(1+i+k)-j}}\Bigg].\nonumber \\
	&& \label{lag4}\\
	\text{Again,} ~ I_2 &=& \int_{a}^{a+b}\exp\Big[-\frac{t^2+\sqrt{2t^2+t^4}}{2}\Big]~dt \nonumber\\
	&=& \sum_{i=0}^{\infty}\frac{(-1)^i}{2^i i!}\sum_{j=0}^{i} {i \choose j}\int_{a}^{a+b}t^{2(i-j)}(2t^2+t^4)^{j/2}~dt \nonumber \\
	&=& \sum_{i=0}^{\infty}\frac{(-1)^i}{2^i i!}\sum_{j=0}^{i} {i \choose j}\int_{a}^{a+b}t^{2(i-j)+j}\cdot2^{j/2}\Big(1+\frac{t^2}{2}\Big)^{j/2}~dt ~;~ \Big|\frac{t^2}{2}\Big| < 1 \nonumber\\
	&=& \sum_{i=0}^{\infty}\frac{(-1)^i}{2^i i!}\sum_{j=0}^{i} {i \choose j}\sum_{k=0}^{\infty} {j/2 \choose k} 2^{j/2 - k}\int_{a}^{a+b}t^{2i-j+2k}~dt \nonumber \\
	&=& \sum_{i=0}^{\infty}\sum_{j=0}^{i}\sum_{k=0}^{\infty}{i \choose j} {j/2 \choose k}\frac{(-1)^i ~2^{j/2 - i - k}}{i!}\Bigg[\frac{(a+b)^{1+2(i+k)-j} - a^{1+2(i+k)-j}}{{1+2(i+k)-j}}\Bigg]. \nonumber\\
	&& \label{lag5}
\end{eqnarray}
 Hence, from eqns. (\ref{lag4}) and (\ref{lag5}), we get\\
$ L.H.S = \frac{e^{-1/2}}{b^2}[I_1 -aI_2] = R.H.S = \hat{\nu_m}$ implies
\begin{align}
	\sum_{i=0}^{\infty}\sum_{j=0}^{i}\sum_{k=0}^{\infty}{i \choose j} {j/2 \choose k}\frac{(-1)^i ~2^{j/2 - i - k}}{i!} \cdot \Bigg[\Bigg(\frac{(a+b)^{2(1+i+k)-j} - a^{2(1+i+k)-j}}{{2(1+i+k)-j}}\Bigg) \nonumber \\ 
	-a\Bigg(\frac{(a+b)^{1+2(i+k)-j} - a^{1+2(i+k)-j}}{{1+2(i+k)-j}}\Bigg)\Bigg] = b^2 e^{1/2}\hat{\nu}_m. \label{lag6} 
\end{align}
\hspace*{0.2in} Approximating up to two terms in (\ref{lag6}), we work out the case for $(i=0,1; k=0)$ to get a second linear equation in $ a $ and $ b $. This is given as: 
\begin{equation}
	3a+2b = -6\sqrt{2}e^{1/2}\hat{\nu}_m. \label{lag7} 
\end{equation}
\hspace*{0.2in} Henceforth, the Lagrange multipliers can be  determined by solving the system of two linear equations in $a$ and $b$, given by (\ref{lag3}) and (\ref{lag7}). This gives us a simplified model for deriving the distribution of $ \hat{\nu}. $ 
\subsection{Validation of model with experimental and field data} 
\begin{figure}[h]
	\includegraphics[width=0.4\textwidth]{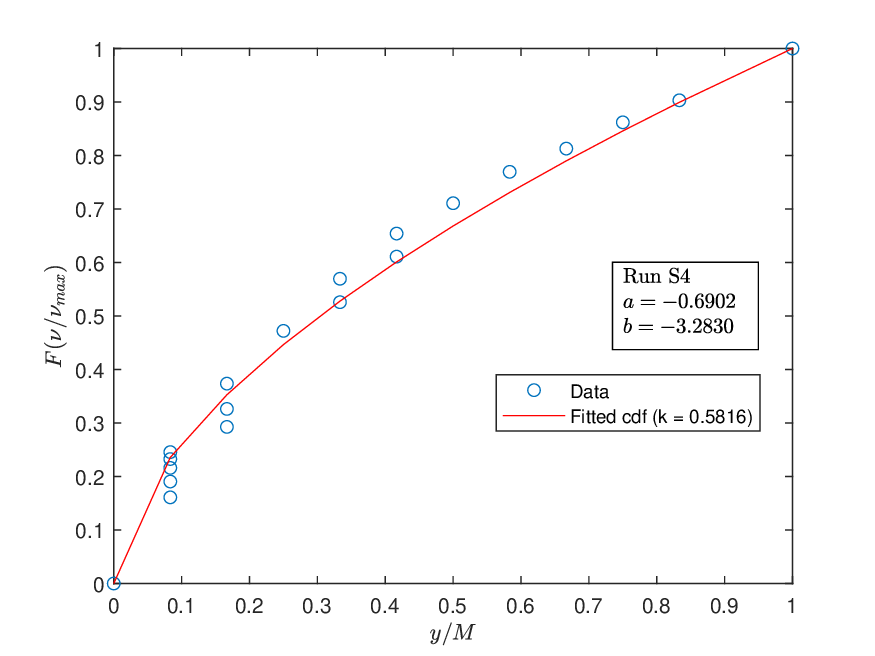} ~~~~~
	\includegraphics[width=0.4\textwidth]{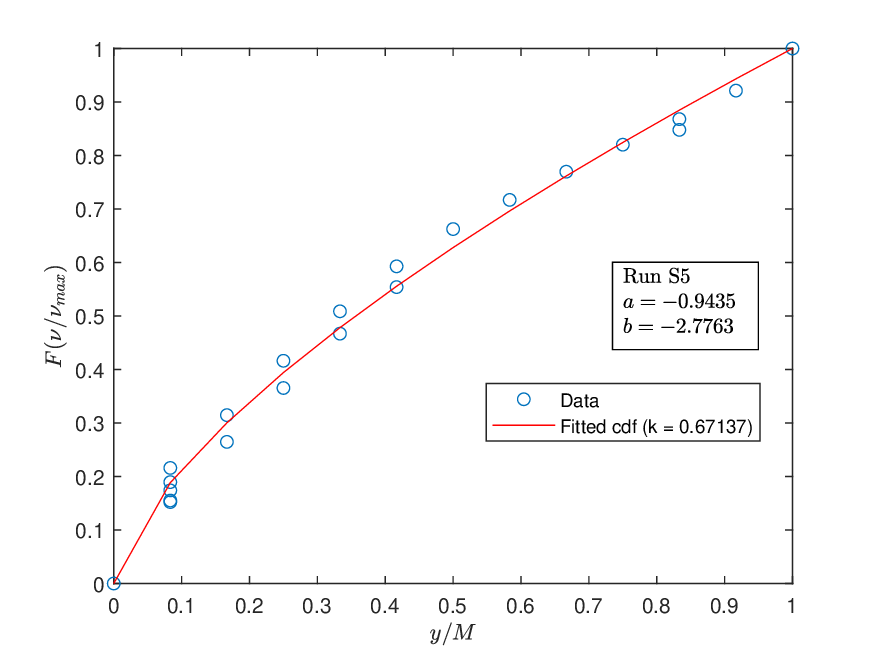}~~~\\
	\includegraphics[width=0.4\textwidth]{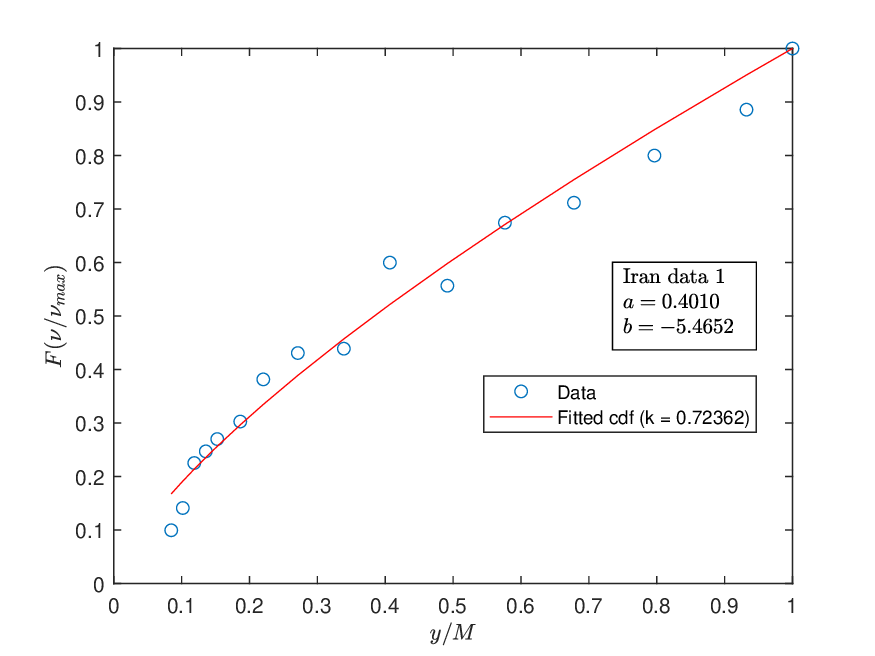}~~~~~ 
	\includegraphics[width=0.4\textwidth]{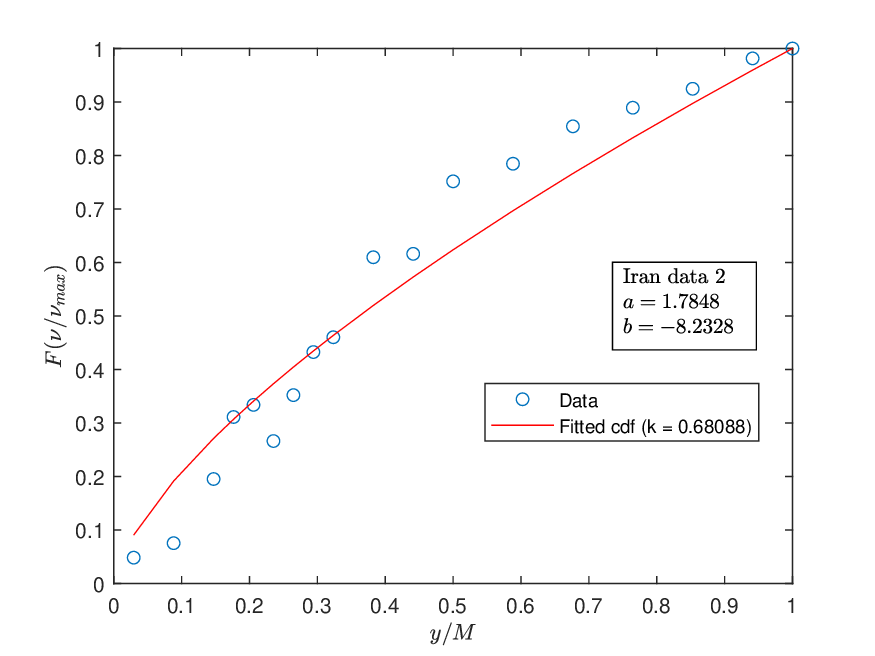} 
	\caption{Validation of the fitted hypothetical cdf (\ref{est2}) with the computed cdf (\ref{est6}) for Einstein \& Chien and Iran series data of Luo }
	\label{cdf_fit}
\end{figure}

\begin{figure}[]
	\includegraphics[width=0.4\textwidth]{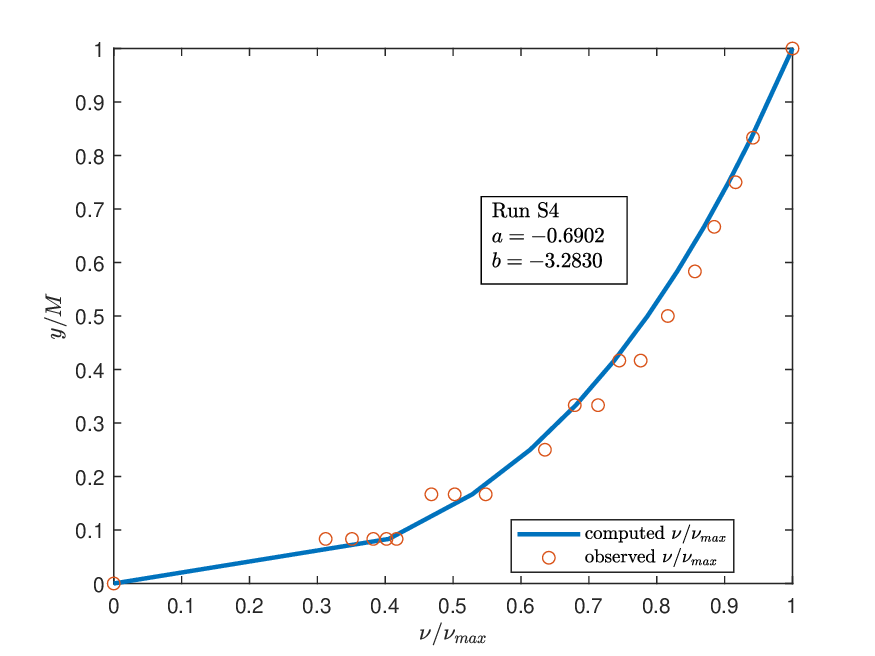} ~~~~~
	\includegraphics[width=0.4\textwidth]{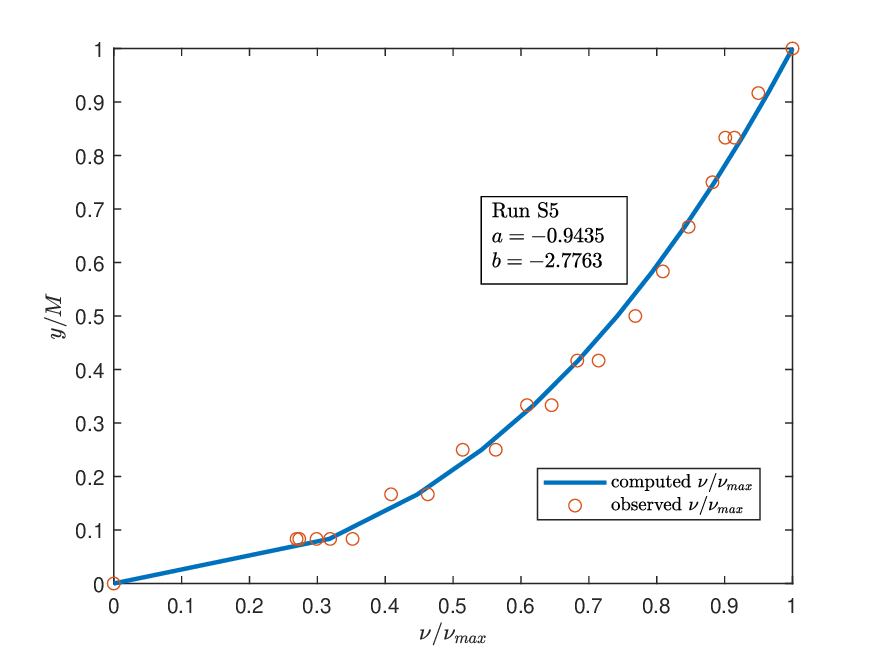}\\
	\includegraphics[width=0.4\textwidth]{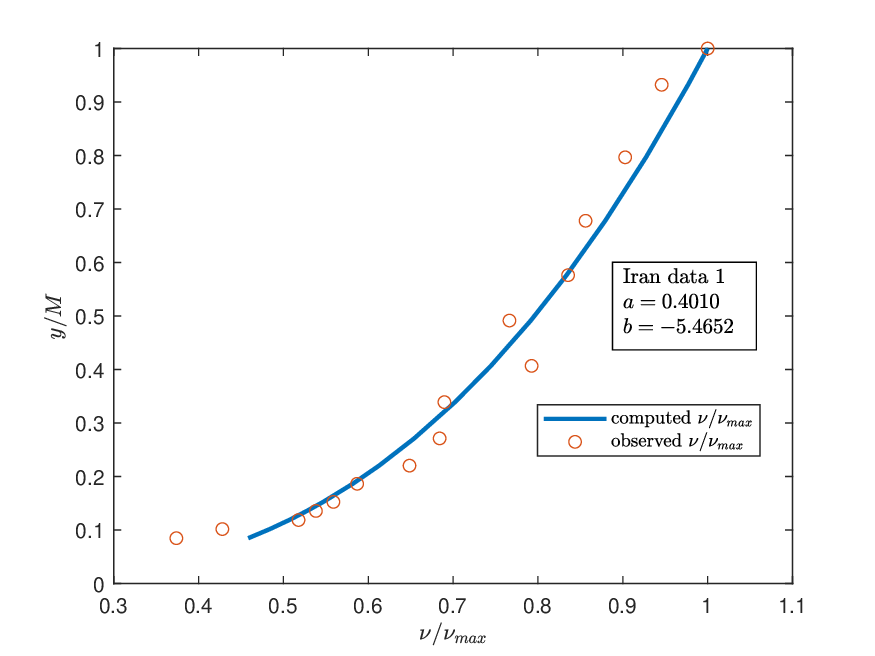} ~~~~~
	\includegraphics[width=0.4\textwidth]{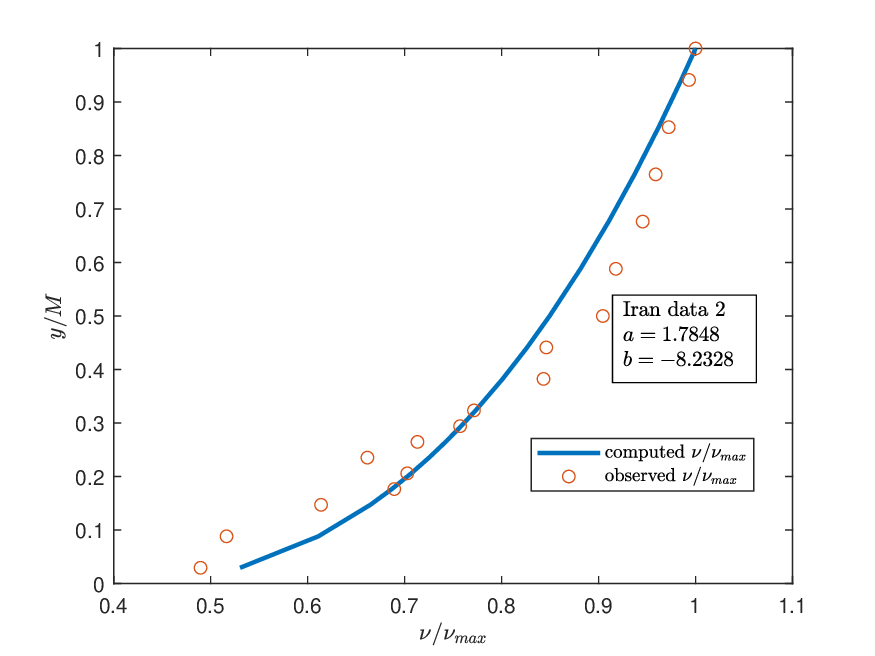} 
	\caption{Normalized velocity ($ \nu/\nu_{max} $) with respect to normalized height ($ y/M $)}
	\label{vel_ht}
\end{figure}

\begin{figure}[h]
	\includegraphics[width=0.4\textwidth]{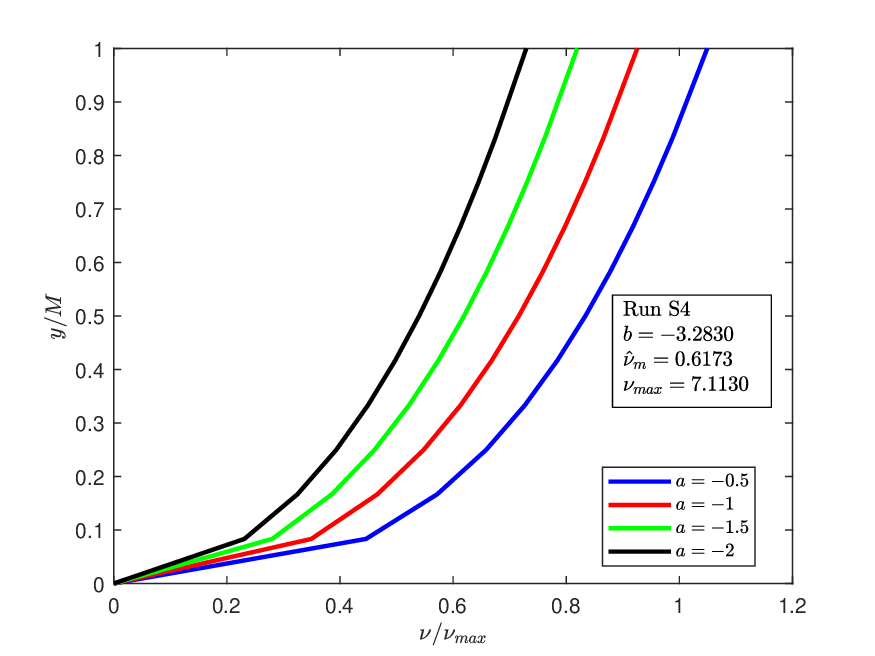} ~~~~~
	\includegraphics[width=0.4\textwidth]{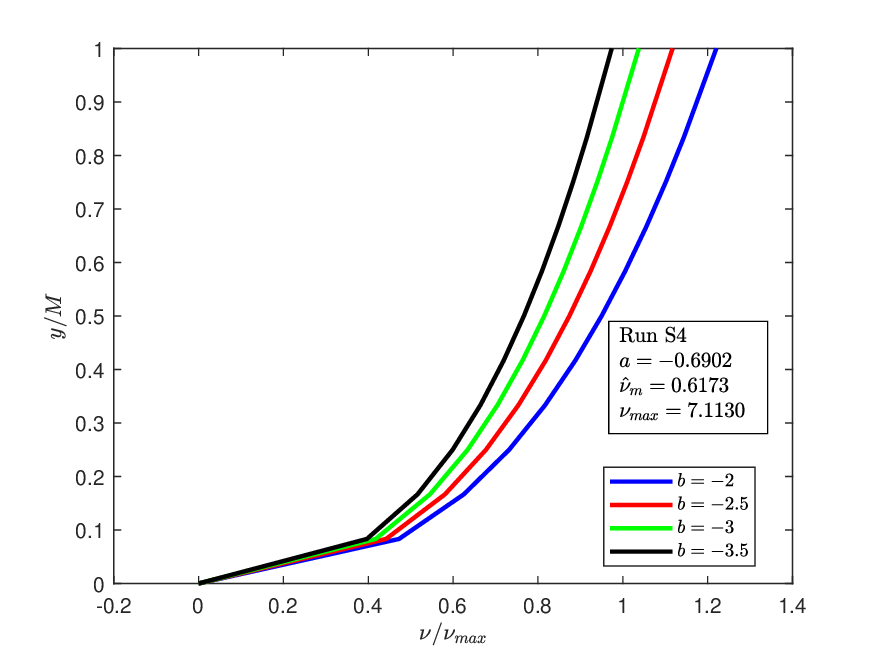}\\
	\includegraphics[width=0.4\textwidth]{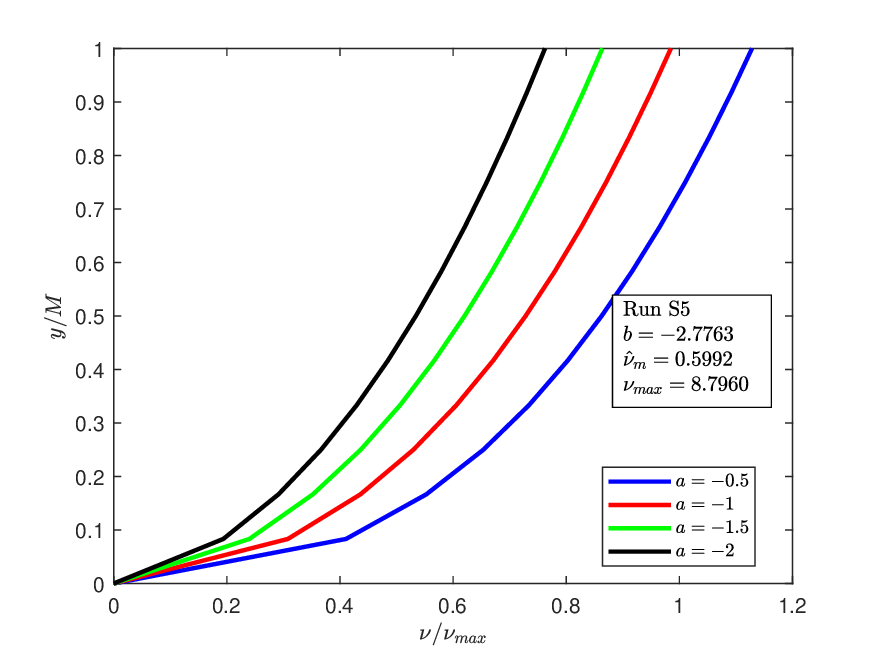} ~~~~~
	\includegraphics[width=0.4\textwidth]{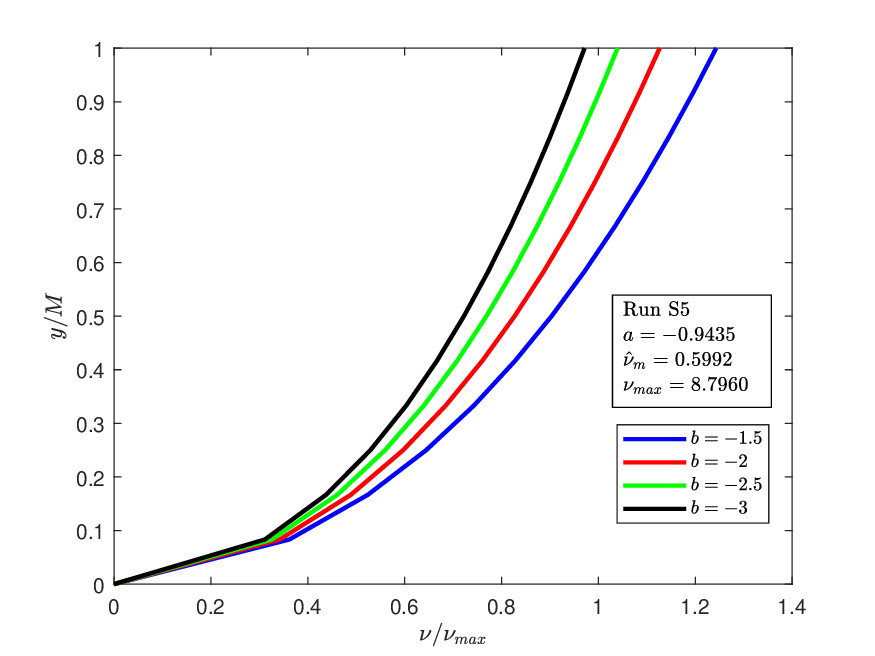}
	\caption{Variation of normalized flow velocity with Lagrange multipliers for Einstein and Chien data}
	\label{lag_1}
\end{figure}

\begin{figure}[]
	\includegraphics[width=0.4\textwidth]{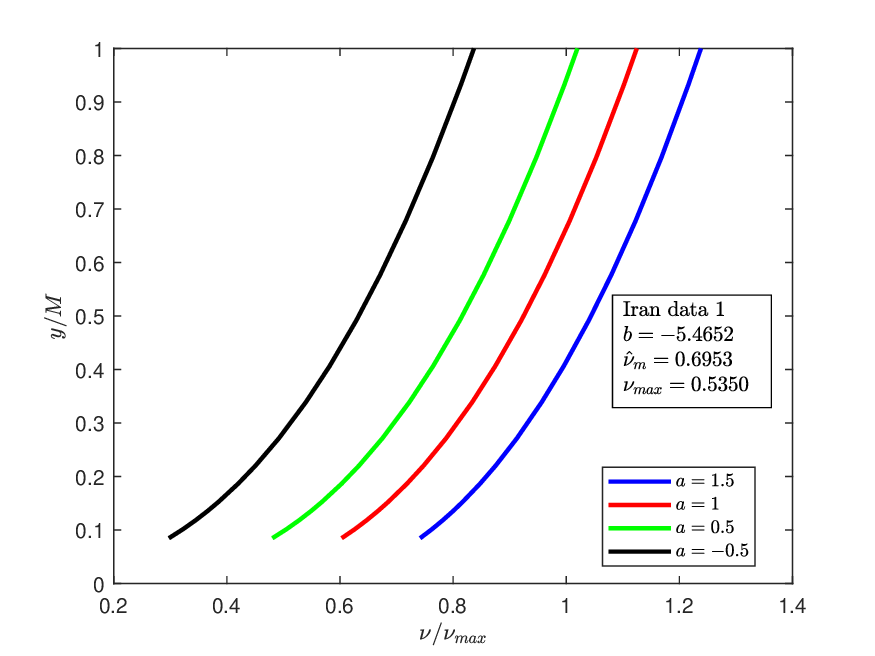} ~~~~~
	\includegraphics[width=0.4\textwidth]{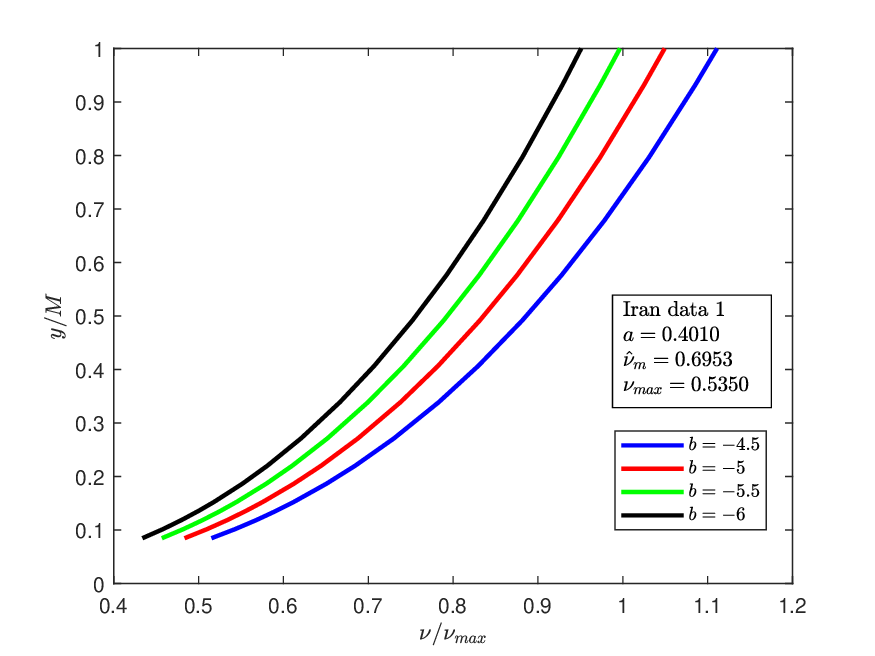}\\
	\includegraphics[width=0.4\textwidth]{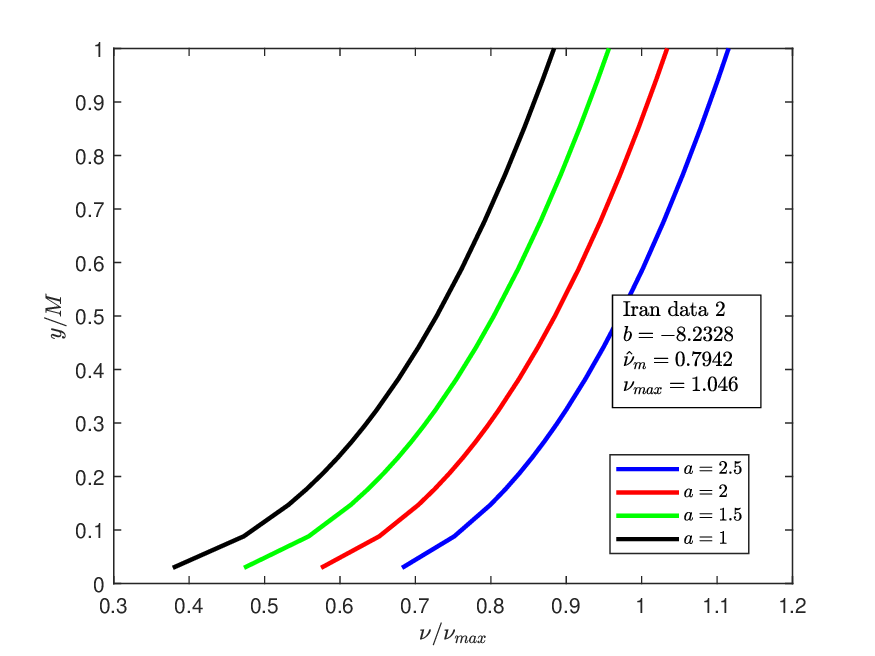} ~~~~~
	\includegraphics[width=0.4\textwidth]{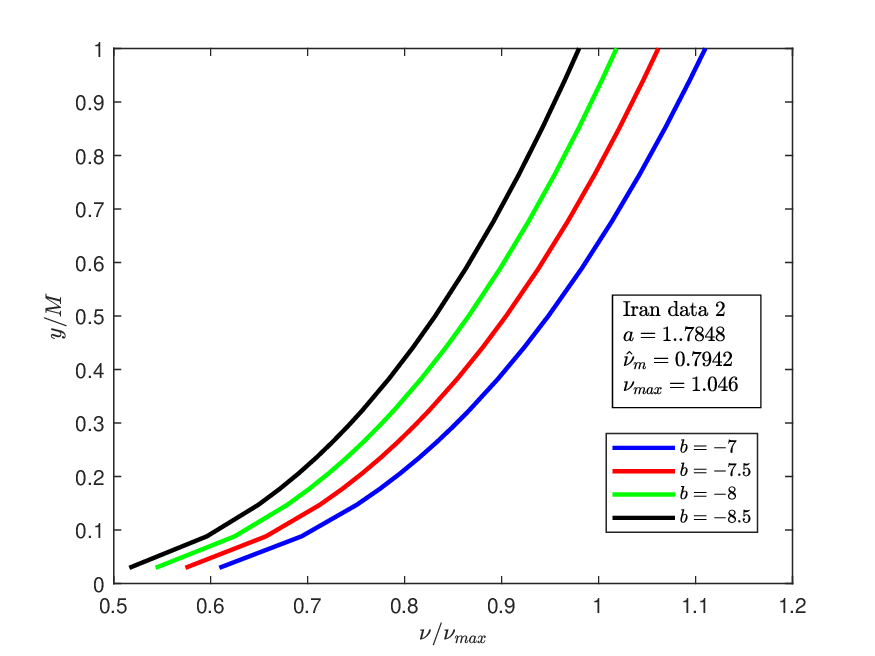}
	\caption{Variation of normalized flow velocity with Lagrange multipliers for Iran series data of Luo}
	\label{lag_2}
\end{figure}

\begin{figure}[h]
	\includegraphics[width=0.4\textwidth]{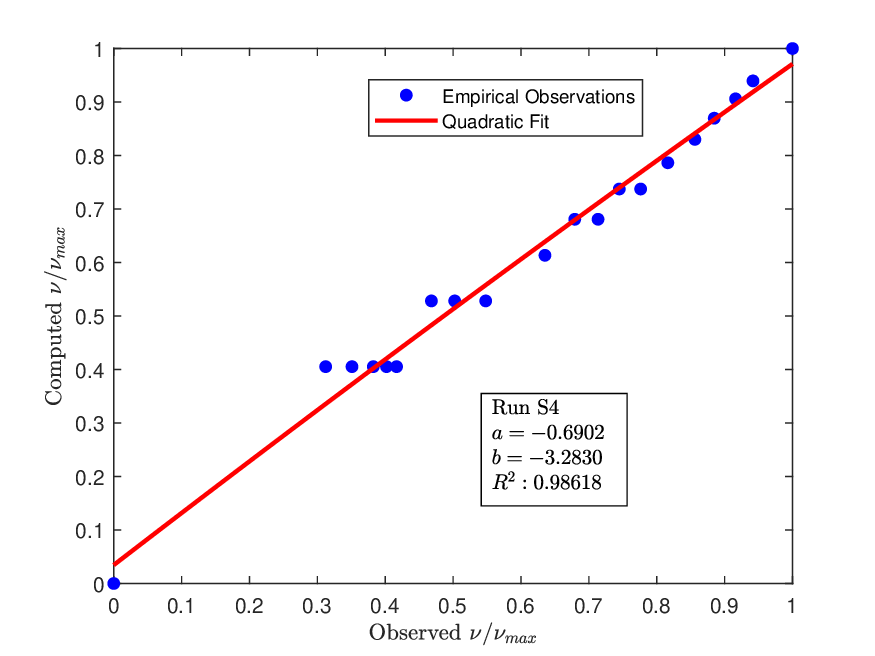} ~~~~~
	\includegraphics[width=0.4\textwidth]{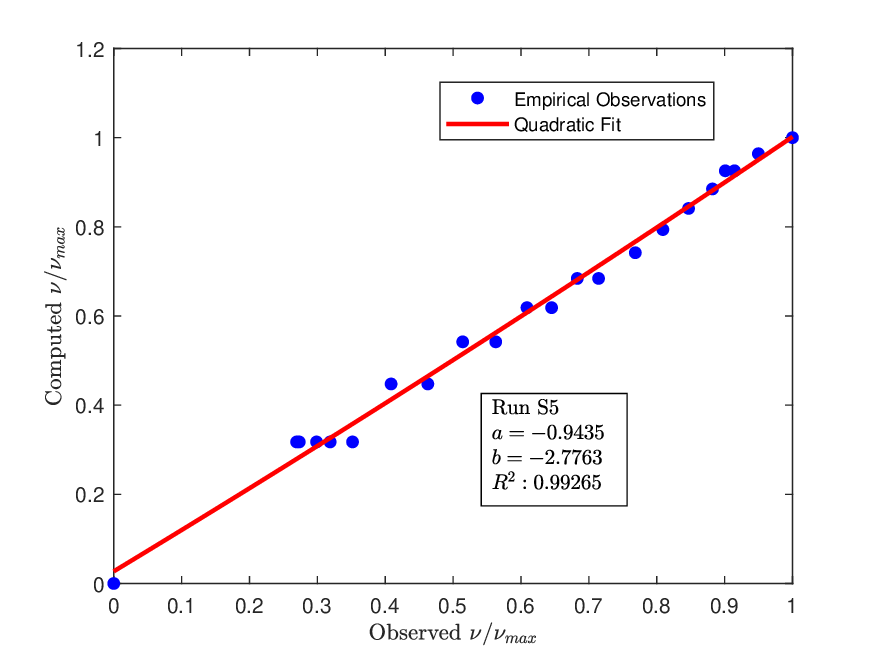} \\
	\includegraphics[width=0.4\textwidth]{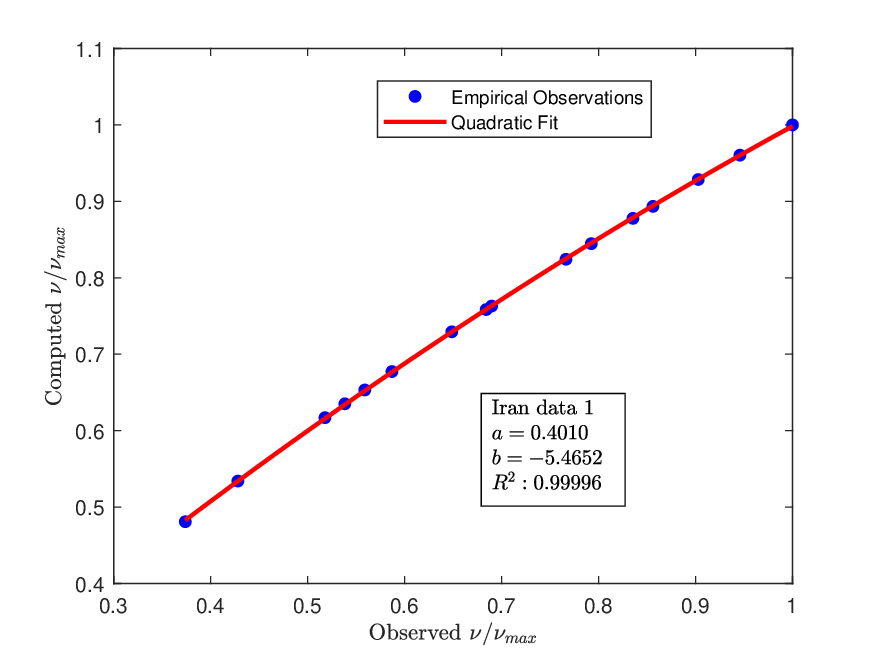} ~~~~~	
	\includegraphics[width=0.4\textwidth]{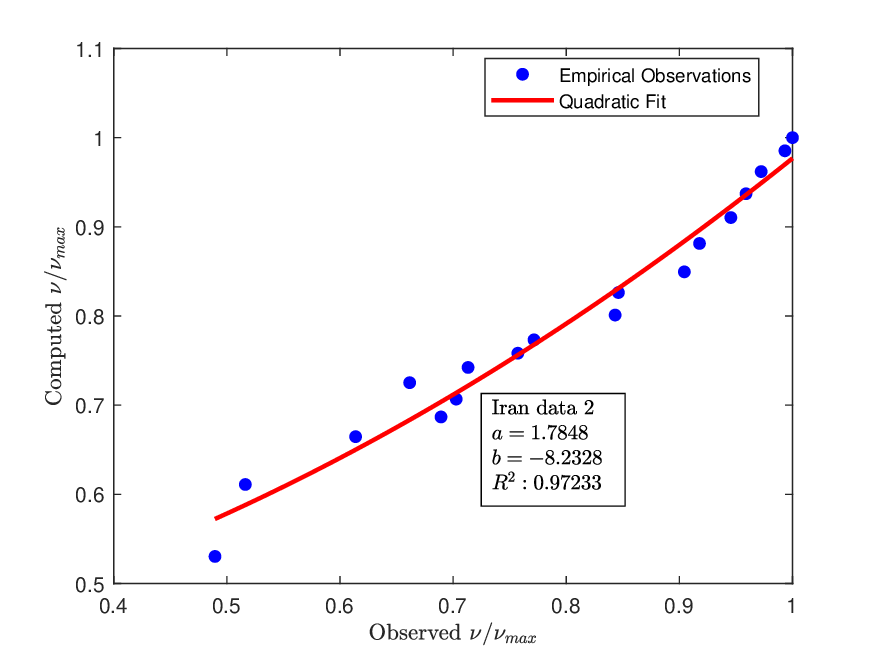} 
	\caption{Regression analysis for computed $ \nu/\nu_{max} $ with respect to observed $ \nu/\nu_{max} $}
	\label{quad_fit}
\end{figure}

The one-dimensional velocity distribution (\ref{est7}) derived from FDE (\ref{2.5}) are investigated using a set of laboratory data collected by Einstein and Chien (1955) and field data collected by Afzalmehr (2008) for channels in Iran. The data sets can be found in Luo (2009). There are 4 sets of velocity data which were used in this study. \\
\hspace*{0.2in} The experimental data collected by Einstein and Chien (1955) were used to analyze the effect of suspended sediment and the roughness of channel bed on the velocity profile near the channel bed. Among a total of 29 runs, 13 were conducted with clear water and the remaining 16 were made with sediment-laden flow. The experiments were conducted in a painted steel flume with a width of 1.006 ft, depth of 1.17 ft and a length of 40 ft. The slope was adjusted with the help of a jack and made to vary within the range 0.0185 to 0.025. The water depth ranged from 0.36 to 0.49 ft and the average velocity of different runs fluctuated from 6.1 ft/s to 8.7 ft/s.
 Data from rivers in Iran (Afzalmehr, 2008) was gathered from wide rectangular channels with clear flow, with observations spanning the entire water depth. Five runs were conducted to evaluate the overall performance of the FDE entropy-based velocity distribution.\\
\hspace*{0.2in} In order to validate the appropriateness and feasibility of the FDE-based velocity distribution model given by (\ref{est7}), firstly we have fitted the FDE-based cdf derived from the pdf of $ \mathcal{\hat{V}} $ represented by (\ref{est6}) with the proposed hypothetical cdf defined by (\ref{est2}) through least square fitting method depicted in Fig. \ref{cdf_fit}. Through this, the value of the fitting parameter $ k $ is determined, which in turn is used to test the proposed velocity distribution model (\ref{est7}) through the selected experimental and field data sets. This is done by plotting the normalized dimensionless velocity and the normalized height $ y/M. $ The results are displayed in Fig. \ref{vel_ht}. From Fig. \ref{vel_ht}, it is evident that the expression for estimation of the most probable velocity distribution among the two resulting cases with positive and negative signs obtained from eqn. (\ref{est7}) should be considered as:
\begin{equation}
	\hat{\nu} = \frac{1}{b}\Bigg[-a-\sqrt{a^2 - 2\sqrt{2}be^{1/2}\left(\frac{y}{M}\right)^k} \Bigg]. \label{val1}
\end{equation}  
\hspace*{0.2in} Next, to check the effect of the Lagrange multipliers with the variation of the normalized velocity $ \hat{\nu}$ and the normalized flow height $ y/M $, the velocity variable $ \hat{\nu} $ is plotted against $ y/M $ for the chosen data sets as shown in Figs. \ref{lag_1} and \ref{lag_2}. Here, we have varied one Lagrange multiplier keeping the other fixed for a given data set. From Figs. \ref{lag_1} and \ref{lag_2}, one can deduce that the vertical velocity variable is more sensitive towards changes in the Lagrange multiplier $ a $ associated with the complete normalization condition. We have that with the same amount of changes in $ a, $ the gradient of velocity variable with respect to height $ y/M $ shows greater fluctuations.\\ 
\hspace*{0.2in} Moreover, to check the correlation between the estimated velocity values obtained from (\ref{est7}) and the observed $ \hat{\nu} $ values obtained from the selected data sets, we have computed the regression coefficients. We have also plotted the correlation between the observed and computed $ \hat{\nu}= \nu/\nu_{max} $ values displayed in Fig. \ref{quad_fit}. The goodness and accuracy of our proposed one-dimensional vertical velocity model represented by eqn. (\ref{est7}) is strongly supported by the regression coefficient $ R^2 $ values which ranges from 0.97233 to 0.99996. Thus, we can find that the $ R^2 $ values obtained for the chosen data sets are very close to 1, thereby showing high correlation between the observed and computed $ \hat{\nu} $ values. This shows that our proposed model has strong prediction accuracy.\\
\hspace*{0.2in} For further validation, we have compared the current FDE-based velocity distribution model (\ref{est7}) with the existing entropy-based models with maximum surface velocity referred to as the Chiu (Shannon entropy), SL (Tsallis entropy), KG (R$\acute{e}$nyi entropy) and KT 1 and KT 2 (fractional Machado entropy) models defined in Chiu (1987), Singh and Luo (2011), Kumbhakar and Ghoshal (2016) and Kumbhakar and Tsai (2023), respectively. For this purpose, we have done error analysis for the normalized velocity $ \hat{\nu} $.
The types of errors considered in this study for
comparison purpose have been defined in Kumbhakar and Tsai (2023). The errors can be described as follows:
\begin{itemize}
	\item Mean relative absolute error (MRAE) = $\frac{1}{n}\sum_{k=1}^{n}\frac{\left|\hat{\nu}_{k(comp)}-\hat{\nu}_{k(obs)}\right|}{\hat{\nu}_{k(obs)}} $
	\item Root-mean-square error (RMSE) = $\mathlarger{\sqrt{\frac{1}{n}\sum_{k=1}^{n}(\hat{\nu}_{k(comp)}-\hat{\nu}_{k(obs)})^2}}$
\end{itemize} 
where $ n $ is the total number of data points in a particular data set, $ \hat{\nu}_{k(comp)} $ and $ \hat{\nu}_{k(obs)} $ denote the computed and observed concentration of sediment at the $ k^{th} $ data point. Since the velocity at the channel bed is taken to be zero, the two errors have been computed by ignoring the data point at the bed region to get finite values. The results are presented in Table \ref{error} (given at the end). \\
\hspace*{0.2in} We show that our model has minimum errors for most of the data sets. The KT model 2 with the momentum coefficient $ \beta $ obtained from eqn. (29) of Kumbhakar and Tsai (2023) also shows better results in many cases as indicated in Table \ref{error}. In conclusion, both our proposed model and KT Model 2 outperform other existing entropy-based models that rely on the maximum surface velocity assumption, demonstrating their superiority. However, KT Model 2 requires the estimation of additional parameters when determining the Lagrange multipliers, making it more complex. In contrast, our model offers a simplified approach by computing the Lagrange multipliers through a system of linear equations. Additionally, our model demonstrates greater efficiency in supporting field data. Therefore, compared to KT Model 2, our proposed model is both simpler and more effective.
\section{Conclusions}
The entropy is used as a quantitative measure of randomness or uncertainty of information related to a real world random phenomena like the hydraulic velocity variable in our case. A plethora of entropy computation methods exist in literature and among them the FDE based on probability density function was chosen. Some univariate distributions popularly used in reliability and survival analysis studies were considered to obtain the closed form expressions of their associated FDEs. However, since closed-form analytical expressions were not feasible for all distributions, some bounds for FDE were explored and its monotonic properties were analyzed. \\
\hspace*{0.2in} Further, the definition of FDE was utilized to formulate a one-dimensional vertical velocity distribution expression for the centre region of wide-open channels, where maximum velocity was assumed to be observed at the surface of the flow. This was done by using the principle of maximum entropy introduced by Jaynes. Through this study, the FDE was found to be quite useful in dealing with the intrinsic robustness of hydrology data by choosing a particular value of the fractional parameter $ \alpha $. A one-parametric cdf was employed to enhance its estimation using a fitting parameter $ k $, which varies based on the characteristics of sediment particles in the flow domain. Additionally, a simplified approach was introduced for determining the Lagrange multipliers in the proposed FDE-based velocity model, improving computational efficiency compared to existing deterministic and entropy-based methods. Finally, the proposed velocity distribution model was validated using both experimental and field data. The low computed error values and high regression coefficients confirmed the model's reliability and effectiveness in representing the velocity distribution of turbulent flows along the vertical profile of wide open channels. Furthermore, a comparison with existing entropy-based models demonstrated the superiority of the proposed FDE-based model, as it yielded the lowest errors in most cases. Notably, it performed particularly well with the field data used in our study.

\section{Disclosure statement}
On behalf of all authors, the corresponding author states that there is no conflict of interest.


\begin{landscape}
	\begin{table}[h]
		\small
		\centering
		\begin{threeparttable}
			\caption{Comparison of the present model with other vertical profile of velocity distribution models}
			\vspace*{0.4in}
			\begin{tabular}{ l c c c c c c}
				\hline
				\textbf{SSC models} &~ \textbf{Proposed model}  &~  \textbf{Chiu model}  &~ \textbf{SL model} &~ \textbf{KG model} &~ \textbf{KT model 1} &~ \textbf {KT model 2}\\
				\hline \hline
				\textbf{\textit{Run S4 of Einstein \& Chien data}} &~ &~ &~ &~ &~ &~ \\ 
				MRAE &~ 0.0529 &~ 0.0291  &~ 0.1046 &~ 0.0264 &~ 0.0880 &~ 0.0239 \\
				RMSE &~ 0.0340 &~ 0.0142  &~ 0.0456 &~ 0.0125 &~ 0.0381 &~ 0.0114 \\
				\textbf{\textit{Run S5 of Einstein \& Chien data}} &~ &~ &~ &~ &~ &~ \\ 
				MRAE &~ 0.0449* &~ 0.0607  &~ 0.1241 &~ 0.0530 &~ 0.1330 &~ 0.0499 \\
				RMSE &~ 0.0243 &~ 0.0265  &~ 0.0460 &~ 0.0240 &~ 0.0495 &~ 0.0234 \\
				\textbf{\textit{Iran series 1 of Luo data}} &~ &~ &~ &~ &~ &~ \\ 
				MRAE &~ 0.0456 &~ 0.0475  &~ 0.0484 &~ 0.0479 &~ 0.0889 &~ 0.0450 \\
				RMSE &~ 0.0336* &~ 0.0381  &~ 0.0422 &~ 0.0393 &~ 0.0601 &~ 0.0363 \\
				\textbf{\textit{Iran series 2 of Luo data}} &~ &~ &~ &~ &~ &~ \\ 
				MRAE &~ 0.0417* &~ 0.0675  &~ 0.0775 &~ 0.0690 &~ 0.0639 &~ 0.0622 \\
				RMSE &~ 0.0384* &~ 0.0641  &~ 0.0658 &~ 0.0634 &~ 0.0588 &~ 0.0584 \\
				\hline			 
			\end{tabular}
			\label{tab3}
			\begin{tablenotes}
				\footnotesize
				\item Note: * indicates the minimum values of each of the errors for a given data source.
			\end{tablenotes}
		\end{threeparttable}
	\end{table}	
	\label{error}
\end{landscape}
\end{document}